\baselineskip=15pt plus 2pt
\magnification =1200
\def\sqr#1#2{{\vcenter{\vbox{\hrule height.#2pt\hbox{\vrule width.#2pt
height#1pt\kern#1pt \vrule width.#2pt}\hrule height.#2pt}}}}
\def\square{\mathchoice\sqr64\sqr64\sqr{2.1}3\sqr{1.5}3}
\centerline{\bf On tau functions for orthogonal polynomials and matrix
models}\par
\vskip.1in
\centerline {Gordon Blower}\par
\centerline {Department of Mathematics and Statistics}\par
\centerline {Lancaster University}\par
\centerline {Lancaster, LA1 4YF}\par
\centerline {England}\par
\centerline {g.blower@lancaster.ac.uk}\par
\vskip.05in
\centerline {13th August 2010}\par
\vskip.05in
\noindent {\bf Abstract.} Let $v$ be a real polynomial of even degree,
and let $\rho$ be the equilibrium probability measure for $v$ with
support $S$; so that, $v(x)\geq 2\int\log \vert x-y\vert \, \rho
(dy)+C_v$ for some constant $C_v$ with equality on $S$. Then $S$ is the
union of finitely many bounded intervals with endpoints $\delta_j$,
and $\rho$ is given by an algebraic weight $w(x)$ on $S$. The
system of orthogonal polynomials for $w$ gives rise to
the Magnus--Schlesinger differential equations.
This paper identifies the $\tau$ function of this system with the
Hankel determinant $\det [\int x^{j+k}\rho (dx)]_{j,k=0}^{n-1}$ of
$\rho$. The solutions of the Magnus--Schlesinger equations are
realised by a linear system, which is used to compute the tau
function in terms of a Gelfand--Levitan equation. The tau function
is associated with a potential $q$ and a scattering problem for the
Schr\"odinger operator with potential $q$. For some
algebro-geometric $q$, the paper solves the scattering
problem in terms of linear systems. The theory extends naturally to
elliptic curves and resolves the case where $S$ has exactly two
intervals.\par
\par
\vskip.05in
\noindent {\bf 1. Introduction}\par
\vskip.05in
\noindent This paper concerns systems of orthogonal polynomials that
arise in random matrix theory, specifically in the theory of the
generalized unitary ensemble [26], and may be described in terms of
electrostatics. We consider a unit of charge to be distributed
along an infinite conducting wire in the presence of an electrical
field. The field is represented by a real polynomial 
$v(x)=\sum_{j=0}^{2N}a_jx^j$ such that $a_{2N}>0$, while the charge is
represented by a Radon probability measure on the real line.\par 
\indent Boutet de Monvel {\sl et al} [5, 29] prove the existence of the
equilibrium distribution $\rho$ that minimises the electrostatic
energy. Under general conditions which include the above $v$, they prove that
there exists a constant $C_v$ such that 
$$v(x)\geq 2\int_S \log \vert x-y\vert\, \rho (dy) +C_v\qquad (x\in {\bf
R})\eqno(1.1)$$
\noindent and that equality holds if and only if $x$ belongs to a
compact set $S$. Furthermore, there exists $g\geq 0$ and 
$$-\infty <\delta_1
<\delta_2\leq \delta_3<\dots <\delta_{2g+2}<\infty \eqno(1.2)$$
\noindent such that 
$$S=\cup_{j=1}^{g+1}[\delta_{2j-1}, \delta_{2j}]\eqno(1.3)$$
\noindent It is a tricky problem, to find $S$ for a given $v$, and [10, Theorem 1.46 and p. 408]
contains some significant results including the bound $g+1\leq
N+1$ on the number of intervals. When $v$ is convex, a relatively
simple argument shows that $g=0$, so there is a single interval [21]. \par

\vskip.05in
\noindent {\bf Definition.} The $n^{th}$ order Hankel determinant for
$\rho$ is
$$D_n=\det\Bigl[ \int_Sx^{j+k}\rho (dx)\Bigr]_{j,k=0}^{n-1}.\eqno(1.4)$$
\vskip.05in
\noindent In section 3 we introduce the system of orthogonal
polynomials for $\rho$, and in  section 4, we regard $D_n$ as a function of $\delta =(\delta_1, \dots ,
\delta_{2g+2})$, and derive a system for differential equations for
$\log D_n$, known as Schlesinger's equations. Let $A(z)$ be a proper rational $2\times 2$ matrix function
with simple poles at $\delta_j$; let $\alpha_j$ be the residue at
$\delta_j$, and suppose that the eigenvalues of $\alpha_j$ are distinct
 modulo the integers for $j=1, \dots , M$. Consider the differential equation
$${{d}\over{dz}}\Phi =A(z)\Phi (z),\eqno(1.5)$$
\noindent and introduce the $1$-form 
$$\Omega (\delta
)={{1}\over{2}}\sum_{j=1}^M{\hbox{trace}}\, {\hbox{residue}}(A(z)^2:
z=\delta_j)\, d\delta_j\eqno(1.6)$$
\noindent to describe its deformations. Then $\Omega$ turns out to be 
closed by [19].\par
\vskip.05in
\noindent {\bf Definition.} (i) The tau function of the deformation 
equations associated with (1.5) is $\tau :{\bf C}^M\setminus
\{{\hbox{diagonals}}\}\rightarrow {\bf C}$ such that 
$d\log \tau =\Omega.$\par
\indent (ii) Given a self-adjoint and trace-class operator $K:L^2 (\rho )\rightarrow
L^2(\rho )$ such that $0\leq K\leq I$, and $P_{(t,\infty )}$ the
orthogonal projection $f\mapsto {\bf I}_{[t, \infty )}f$, the tau
function of $K$ is $\tau (t)=\det (I-P_{[t, \infty )}K)$. (The
definitions in (i) and (ii) are reconciled by Proposition 3.3.)\par  
\vskip.05in
\indent In section 5 use the results from preliminary
sections to prove that $D_n$ gives the
appropriate $\tau$ function for Schlesinger's equations. As an illustration which is of importance
in random matrix theory, we calculate the tau function explicitly when $\rho$ is the
semicircular law. When $S$ is the union of two intervals, the
Schlesinger equations reduce to the Painlev\'e VI equation, as we
discuss on section 7. Okamoto derived $\tau$ functions for other
Painlev\'e equations in [27]. See also [15] and [2].\par
\indent In sections 2,3 and 4, we develop standard arguments, then our analysis follows that of Chen and Its [8], who considered the $\rho$ that is
analogous to the Chebyshev distribution on multiple intervals. 
 Chen and Its found their tau function explicitly in
terms of theta functions on a hyperelliptic Riemann surface; in this
paper, we show by means of scattering theory why 
such solutions emerge.\par
\indent In terms of definition (ii), we have a real $\tau (x)$ which is associated with a compactly supported 
real potential $q(x)=-2{{d^2}\over{dx^2}}\log \tau (2x)$ and hence the
Schr\"odinger differential operator $-{{d^2}\over{dx^2}}+q(x)$. One associates
with each smooth $q$ a scattering
function $\phi$; then one analyses the spectral data of 
$-{{d^2}\over{dx^2}}+q(x)$ in
terms of $\phi$, with a view
to recovering $q$. The Gelfand--Levitan integral equation links
$\phi$
with $q$.\par

\indent In random matrix theory, tau functions are introduced
alongside integrable kernels that describe the distribution of eigenvalues
of random matrices, especially the generalized unitary ensemble; 
see [14, 26, 31]. Let $X$ be a $n\times n$ complex Hermitian matrix, let
$\lambda =(\lambda_1\leq \lambda_2\leq \dots \leq \lambda_n)$ be the
corresponding eigenvalues, listed according to multiplicity, and 
consider the potential $V(X)=n^{-1}\sum_{j=1}^n v(\lambda_j)$.
Now let $dX$ be the product of Lebesgue measure on the entries that are
on or above the leading diagonal of $X$; then there exists
$0<Z_n<\infty$ such that 
$$\nu_n^{(2)}(dX)=Z_n^{-1}\exp\bigl( -n^2V(X)\bigr) dX\eqno(1.7)$$
\noindent defines a probability measure on the $n\times n$ complex
Hermitian matrices. There is a natural action of the unitary
group $U(n)$ on $M_n$ given
by $(U,X)\mapsto UXU^\dagger$, which leaves $\nu_n^{(2)}$ invariant.
Hence $\nu_n^{(2)}$ is the generalized unitary ensemble with potential
$v$.\par
\indent There exists a constant $\zeta_n$ such that
$\rho_n(dx)=\zeta_n^{-1}e^{-nv(x)}\, dx$ defines a probability
measure on ${\bf R}$; then we let $E_k^{\rho_n}$ be the orthogonal
projection onto ${\hbox{span}}\{ x^j: j=0, \dots ,k-1\}$ in
$L^2(\rho_n)$. The eigenvalue distribution satisfies
$$\int{{1}\over{n}}\sharp \{j: \lambda_j(X)\leq t\}\nu_n^{(2)}(dX)=\det
(I-{\bf I}_{(t, \infty )}E_{n}^{\rho_n})\eqno(1.8)$$
\noindent where the right-hand side can be expressed in terms of
Hankel determinants. For large $n$, most
of the eigenvalues actually lie in $S$ by results of [5, 26]. Moreover, there
exists a trace-class operator $K$ on $L^2(\rho )$ such that $0\leq
K\leq I$ and 
$$\det (I-{\bf I}_{(t, \infty )}E_{n}^{\rho_n}))\rightarrow 
\det (I-{\bf I}_{(t, \infty )}K)\qquad (n\rightarrow\infty
);\eqno(1.9)$$ 
\noindent call this limit $\tau (t)$. Tracy and Widom [31] showed how to express such
determinants in terms of systems of differential equations and Hankel
operators.\par
\indent We introduce the matrix
$$J=\left[\matrix{0&-1\cr 1&0\cr}\right],$$
\noindent and apply a simple gauge transformation to (1.5). Then for a sequence of real symmetric $2\times 2$ matrices
 $J\beta_k(n)$, we consider solutions of the differential equation
$$J{{dZ}\over{dx}}=\sum_{k=1}^{2g+2}
{{J\beta_k(n)}\over{x-\delta_k}}Z,\eqno(1.10)$$
$$Z(x)\rightarrow 0\qquad (x\rightarrow \delta_j),$$
\noindent and form the kernel
$$K(x,y)={{Z(y)^\dagger JZ(x)}\over{y-x}}.\eqno(1.11)$$
\noindent We show that the properties of $K$ depend crucially
upon the sequence of signatures of the matrices 
$(\delta_j-\delta_k)J\beta_k(n)$. In Theorem 8.3, we introduce a symbol function $\phi$ from $Z$, a constant
signature matrix $\sigma$ and a Hankel operator $\Gamma_\phi$ such that
$K=\Gamma_\phi^\dagger\sigma\Gamma_\phi.$\par 

\indent In section 9 we introduce $\phi$ from (1.11), express $\phi$ in
terms of a linear system as in [4] and
hence obtain a matrix Hamiltonian $H(x)$ such
that
$$\tau (2x)=\exp\Bigl( -\int_x^\infty {\hbox{trace}}\, H(u)\,
du\Bigr),\eqno(1.12)$$
\noindent and prove that $q(x)$ is
meromorphic on a region. We regard $-{{d^2}\over{dx^2}}+q$ as integrable
if $-f''+qf=\lambda f$ can be solved by quadratures for typical
$\lambda$. This imposes severe restrictions upon $q$; indeed, 
Gelfand, Dikij and Its [12, 6] showed that the integrable cases arise from finite-dimensional Hamiltonian
systems. In sections 10, 11 and 12 we consider cases in which $-f''+qf=\lambda f$ has a meromorphic
general solution for all $\lambda$, and $q$ satisfies one of the
following conditions:\par
\indent (i) $q$ is rational and bounded at infinity;\par
\indent (ii) $q$ is of
rational character on an elliptic curve;\par
\indent (iii) $q$ is the restriction of an abelian function to a straight line in the Jacobian
of a hyperelliptic Riemann surface.\par
\indent In cases (ii) and (iii), the corresponding Schr\"odinger
equation has a spectrum with only finitely many gaps. In (i) and (ii), we introduce a linear system $(-A,B,C)$ so as to 
realise $\phi (x)=Ce^{-xA}B,$ and use the operators $A,B,C$ to solve the Gelfand--Levitan equation. Thus
we obtain explicit expressions for $\phi$ and $\tau$. In (iii), we can
do likewise under further hypotheses.\par
\vskip.05in
\noindent {\bf 2. The equilibrium measure}\par
\vskip.05in
\indent Given the special form of the potential, the equilibrium
measure and its support satisfy special properties. To describe these,
we introduce the polynomial $u$ of degree $2N-2$ by
$$u(z)=\int_S {{v'(z)-v'(x)}\over{z-x}}\rho (dx)\eqno(2.1)$$
\noindent and the Cauchy transform of $\rho$ by
$$R(z)=\int_S {{\rho (dx)}\over{x-z}}\qquad (z\in {\bf C}\setminus
S)\eqno(2.2)$$
\noindent and the weight
$$w(x)=2Na_{2N}\Bigl({-Q(x)
\prod_{j=1}^{2g+2}(x-\delta_{2j-1})(x-\delta_{2j})}\Bigr)^{1/2}\eqno(2.3)$$
\noindent where $Q(x)$ is a product of monic irreducible quadratic
factors such that $w(x)^2=4u(x)-v'(x)^2.$\par
\vskip.05in
\noindent {\bf Proposition 2.1} {\sl (i) The Cauchy transform is the
algebraic function that satisfies
$$R(z)^2+v'(z)R(z)+u(z)=0\eqno(2.4)$$
\noindent and $R(z)\rightarrow 0$ as $z\rightarrow\infty$. There
exist nonzero polynomials $u_0,u_1,u_2$ such that $u_0R'=u_1R+u_2.$\par
\indent (ii) The support of $\rho$ is 
$$S=\{ x\in {\bf R}: 4u(x)-v'(x)^2\geq 0\}.\eqno(2.5)$$
\indent (iii) $\rho$ is absolutely continuous and
the Radon--Nikodym derivative satisfies}
$${{d\rho}\over{dx}}={{1}\over{2\pi }}{\bf I}_{S}(x)w(x)\eqno(2.6)$$
\noindent {\sl where $2\pi =\int_Sw(t)dt$ and $w(x)\rightarrow 0$
as $x$ tends to an endpoint of $S$.}\par
\vskip.05in
\noindent {\bf Proof.} (i) The quadratic equation is due to Bessis, 
Itzykson and Zuber, and is proved in the required form in [28]. One can easily
deduce that $R$ satisfies a first-order linear differential equation
with polynomial coefficients.\par
\indent (ii) Pastur [28] shows that the support is those real $x$ such that
$$\vert
v'(x)+\sqrt{v'(x)^2-4u(x)}\vert^2 =4u(x),\eqno(2.7)$$
\noindent and this condition reduces to $4u(x)\geq v'(x)^2$ and
$u(x)\geq 0$, where the former inequality implies the latter. The polynomial
$4u(x)-v'(x)^2$ has real zeros $\delta_1, \dots , \delta_{2g+2}$, and
may additionally have pairs of complex conjugate roots, which we list
as $\delta_{2g+3}, \dots ,\delta_{4N-2}$ with regard to multiplicity. Hence we can
introduce $w$ as above such that $4u(x)-v'(x)^2=w(x)^2.$\par 
\indent (iii) From (i) we deduce that 
$$R(\lambda )={{1}\over{2\pi i}}\int_S
{{\sqrt{4u(t)-v'(t)^2}}\over{t-\lambda }}dt\eqno(2.8)$$
\noindent since both sides are holomorphic on ${\bf C}\setminus S$,
vanish at infinity and have the same jump across $S$. By Plemelj's
formula, we deduce that
$$v'(\lambda )=2{\hbox{p.v.}}\int_S
{{\sqrt{4u(t)-v'(t)^2}}\over{\lambda
-t}}{{dt}\over{2\pi}}\qquad (\lambda\in S).\eqno(2.9)$$
See [28, 5, 29]. This gives the required expression for $\rho$.\par
\rightline{$\square$}\par
\vskip.05in
\noindent {\bf 3. Orthogonal polynomials}\par
\vskip.05in
\indent First we introduce orthogonal polynomials for $\rho$, then
 the corresponding differential equations. Let $(p_j)_{j=0}^\infty$ 
be the sequence of monic orthogonal
polynomials in $L^2 (\rho )$, where $p_j$ has degree $j$ and let $h_j$ be the constants
such that
$$\int_S p_j(x)p_k(x)\rho (dx)=h_j\delta_{jk};\eqno(3.1)$$
and let $(q_j)_{j=1}^\infty$ be the monic polynomials of the second
kind, where
$$q_j(z)=\int_S{{p_j(z)-p_j(x)}\over{z-x}}\rho (dx)\eqno(3.2)$$
\noindent has degree $j-1$.  On account of Proposition 2.1, the
orthogonal polynomials are semi classical in Magnus's sense [22], although the weight
typically lives on several intervals. The following result is standard in the
theory of orthogonal polynomials; see [8].\par
\vskip.05in
\noindent {\bf Lemma 3.1} {\sl Let $c_n=h_n/h_{n-1}$ and
$b_n=h_n^{-1}\int_S xp_n(x)^2\rho (dx).$ Then\par
\indent (i) the polynomials $(p_n)_{n=0}^\infty$ satisfy the
three-term recurrence
relation
$$xp_n(x)=p_{n+1}(x)+b_{n+1}p_n(x)+c_np_{n-1}(x);\eqno(3.3)$$
\indent (ii) the polynomials $(q_j)_{j=1}^\infty$ likewise satisfy
(3.3);\par
\indent (iii) the Hankel determinant of (1.4) satisfies}
$$D_n=h_0h_1\dots h_{n-1}.\eqno(3.4)$$
\vskip.05in
\indent We introduce also 
$$Y_n(z)=\left[\matrix{ p_n(z)& \int_S{{p_n(t)\rho (dt)}\over{z-t}}\cr
{{p_{n-1}(z)}\over{h_{n-1}}}& {{1}\over{h_{n-1}}}\int_S {{p_{n-1}(t)\rho
(dt)}\over {z-t}}\cr}\right]\eqno(3.5)$$
\noindent and 
$$V_n(z)=\left[\matrix{ z-b_{n+1}&-h_n\cr
1/h_n&0\cr}\right].\eqno(3.6)$$
\vskip.05in
\noindent {\bf Proposition 3.2} {\sl (i) The matrices satisfy the
recurrence relation
$$Y_{n+1}(z)=V_n(z)Y_n(z).\eqno(3.7)$$ 
\indent (ii) The matrix $Y_n(z)$ is invertible,
and $\det Y_n(z)=1$.} 
\vskip.05in
\noindent {\bf Proof.} (i) This follows from (i) and (ii) of the
Lemma 3.1.\par
\indent (ii) This follows by induction, where the induction step
follows from the recurrence relation in (i).\par
\rightline{$\square$}\par
\vskip.05in
\indent We restrict $\rho$ restrict to $(-\infty
,t)\cap S$ and let
$$\mu_{j}(t)=\int_{S\cap (-\infty ,t)}x^{j}\rho (dx)\eqno(3.8)$$
\noindent be the $j^{th}$ moment; the corresponding Hankel determinant is
$$D_{n+1}(t)=\det \bigl[\mu_{j+k}(t)\bigr]_{j,k=0}^{n}.\eqno(3.9)$$
\noindent  Let $E_n:L^2(\rho )\rightarrow {\hbox{span}}\{ x^k:k=0,
\dots ,n-1\}$ be the orthogonal projection; we also introduce the
projection $P_{(t, \infty )}$ on $L^2(\rho )$ given by multiplication
$f\mapsto {\bf I}_{(t, \infty )} f$, where ${\bf I}_{(t, \infty )}$
denotes the indicator function of $(t, \infty )$. \par
\vskip.05in
\noindent {\bf Proposition 3.3} {\sl The tau function of
$E_{n+1}$ satisfies}
$$\det (I-E_{n+1}P_{(t, \infty
)})={{D_{n+1}(t)}\over{D_{n+1}}}.\eqno(3.10)$$
\vskip.05in
\noindent {\bf Proof.} We introduce an upper triangular matrix
$[a_{\ell, j}]_{j, \ell =0}^n$ with ones on the leading diagonal 
such that $p_j(x)=\sum_{\ell =0}^n a_{\ell j}x^\ell$. Then we can compute
$$\eqalignno{\det [\mu_{j+k}(t)]_{j,k=0}^n &=\det [a_{\ell
,j}]_{\ell ,j=0}^n\det [\mu_{j+k}(t)]_{j,k=0}^n\det [a_{k,
m}]_{k,m=0}^m\cr
&=\det\Bigl[ \int_{-\infty}^t p_j(x)p_k(x)\rho
(dx)\Bigr]_{j,k=0}^n&(3.11)\cr}$$
\indent We can also express the operators on $L^2(\rho )$ as matrices
with respect to the orthonormal basis $(p_j/\sqrt{h_j})_{j=0}^n$, and
we find
$$E_{n+1}-E_{n+1}P_{(t, \infty )}E_{n+1}\leftrightarrow 
\Bigl[{{1}\over{\sqrt{h_jh_k}}}\int_{-\infty }^t p_j(z)p_k(z)\rho
(dz)\Bigr]_{j,k=0}^n\eqno(3.12)$$
\noindent so that 
$$\det\Bigl[ \int_{-\infty}^t p_j(x)p_k(x)\rho (dx)\Bigr]_{j,k=0}^n=\det
(E_{n+1}-E_{n+1}P_{(t, \infty )}E_{n+1})h_0\dots h_n.\eqno(3.13)$$
\noindent We deduce that 
$$\det [\mu_{j+k}(t)]_{j,k=0}^n=\det (E_{n+1}-E_{n+1}P_{(t, \infty
)}E_{n+1})D_{n+1}.\eqno(3.14)$$
\rightline{$\square$}\par
\vskip.05in

\noindent {\bf 4. Schlesinger's equations and recurrence relations}\par
\vskip.05in
\noindent Invoking Proposition 3.2(ii), we introduce the matrix function
$$A_n(z)=Y'_n(z)Y_n(z)^{-1}+Y_n(z)\left[\matrix{ 0&0\cr 0&
-w'(z)/w(z)\cr}\right] Y_n(z)^{-1}.\eqno(4.1)$$
\noindent The basic properties of $A_n(z)$ are stated in (i) of the
following Lemma, while (ii) gives detailed information that we need in
the subsequent proof of Theorem 5.1.\par 
\vskip.05in
\noindent {\bf Lemma 4.1} {\sl (i) Let $v'(z)^2-4u(z)$ have zeros
at $\delta_j$ for $j=1, \dots , 4N-2$. Then $A_n(z)$ is a proper rational function so
that
$$A_n(z)=\sum_{j=1}^{4N-2} {{\alpha_j(n)}\over{z-\delta_j}},\eqno(4.2)$$
\noindent where the residue matrices $\alpha_j(n)$ depend implicitly upon $\delta.$\par
 \indent (ii) The $(1,2)$ and diagonal entries of the residue matrices 
satisfy} 
$$\eqalignno{\sum_{k=1}^{4N-2} \alpha_k(n)_{12}&=0;&(4.3)\cr
\sum_{k=1}^{4N-2}\bigl(
\alpha_k(n)_{11}-\alpha_k(n)_{22}\bigr)&=2(n+N)-1;&(4.4)\cr 
\sum_{k=1}^{4N-2}\delta_k\alpha_k(n)_{12}&=-2h_n(n+N).&(4.5)\cr}$$
\vskip.05in
\noindent {\bf Proof.} (i) The defining equation (4.1) for $A_n(z)$ may be
written more explicitly as  
$$\left[\matrix{ p_n'(z)& -\int_S{{p_n(t)w(t)dt}\over{(z-t)^2}}\cr
{{p_{n-1}'(z)}\over{h_{n-1}}}& -{{1}\over{h_{n-1}}}\int_S
{{p_{n-1}(t)w(t)dt}\over {(z-t)^2}}\cr}\right]$$
$$=A_n(z)\left[\matrix{ p_n(z)&\int_S{{p_n(t)w(t)dt}\over{z-t}}\cr
{{p_{n-1}(z)}\over{h_{n-1}}}& {{1}\over{h_{n-1}}}\int_S
{{p_{n-1}(t)w(t)dt}\over{ z-t}}\cr}\right] +\left[\matrix{ 0&
{{w'(z)}\over{w(z)}}\int_S
{{p_{n-1}(t)w(t)dt}\over{z-t}}\cr
0&{{w'(z)}\over{h_{n-1}w(z)}}\int_S
{{p_{n-1}(t)w(t)dt}\over{z-t}}\cr}\right].\eqno(4.6)$$
\noindent By considering the entries, we see that $A_n(z)$ is a proper
rational function with possible simple poles at the $\delta_j$, as in
(4.2). Hence
we have a Laurent expansion
$$A_n(z)={{1}\over{z}}\sum_{k=1}^{4N-2}
\alpha_k(n)+{{1}\over{z^2}}\sum_{k=1}^{4N-2}\delta_k\alpha_k(n)+O\Bigl(
{{1}\over{z^3}}\Bigr)\qquad (z\rightarrow\infty ).\eqno(4.7)$$
\indent (ii) First we
compute the $(1,2)$ entry of $A_n(z)$, namely 
$$\eqalignno{A_n(z)_{12}&=-p_n'(z)\int_S{{p_n(t)w(t)dt}\over{z-t}}-p_n(z)\int_S
{{p_n(t)w(t)dt}\over
{(z-t)^2}}-{{w'(z)}\over{w(z)}}p_n(z)\int_S{{p_n(t)w(t)dt}\over
{z-t}}\cr
&=-{{n}\over{z^2}}\int_S
t^np_n(t)w(t)dt-{{(n+1)p_n(z)}\over{z^{n+2}}}\int_St^np_n(t)w(t)\,
dt\cr
&\qquad -{{w'(z)}\over{w(z)}}{{p_n(z)}\over{z^{n+1}}}\int_St^np_n(t)
w(t)dt+O\Bigl(
{{1}\over{z^3}}\Bigr)&(4.8)\cr}$$
\noindent and we can reduce these terms to
$$A_n(z)_{12}={{-h_{n}n}\over{z^2}}-{{h_n(n+1)}\over{z^2}}-{{h_n(2N-1)}\over{z^2}}+O\Bigl(
 {{1}\over{z^3}}\Bigr),\eqno(4.9)$$
\noindent which gives (4.3) and (4.5).\par

\indent Next, the $(2,2)$ entry of $A_n(z)$ is
$$\eqalignno{A_n(z)_{22}&=-{{p_{n-1}'(z)}\over{h_{n-1}}}\int_S
{{p_n(t)w(t)dt}\over{z-t}}-{{p_n(z)}\over{h_{n-1}}}\int_S{{p_{n-1}(t)w(t)dt}\over{(z-t)^2}}\cr
&\qquad
-{{w'(z)}\over{w(z)}}{{p_n(z)}\over{h_{n-1}}}\int_S{{p_{n-1}(t)w(t)dt}
\over{z-t}}\cr
&=-{{(n-1)z^{n-2}}\over{h_{n-1}z^{n+1}}}\int_St^np_n(t)w(t) dt-{{p_n(z)}
\over{h_{n-1}z^{n+1}}}\int_S nt^{n-1}p_{n-1}(t)w(t)dt\cr
&\qquad -{{2N-1}\over{z}}{{p_n(z)}\over{z^n}}{{1}\over{h_{n-1}}}\int_S
t^{n-1}p_{n-1}(t)w(t)dt+O\Bigl({{1}\over{z^2}}\Bigr)\cr
&={{1-n-2N}\over{z}}+O\Bigl({{1}\over{z^2}}\Bigr)&(4.10)\cr}$$
\indent Similarly, the $(1,1)$ entry is
$$\eqalignno{A_n(z)_{11}&={{p_n'(z)}\over{h_{n-1}}}\int_S
{{p_{n-1}(t)w(t)dt}\over{z-t}}+{{p_{n-1}(z)}\over{h_{n-1}}}
\int_S{{p_n(t)w(t)dt}\over{(z-t)^2}}\cr
&\quad +{{w'(z)}\over{w(z)}}{{p_{n-1}(z)}\over{h_{n-1}}}\int_S
{{p_n(t)w(t)dt}\over{z-t}}\cr
&={{p_n'(z)}\over{h_{n-1}z^n}}\int_St^{n-1}p_{n-1}(t)w(t)dt+
{{(n+1)p_{n-1}(z)}\over{h_{n-1}z^{n+2}}}\int_S t^np_n(t)w(t)dt\cr
&\quad +{{w'(z)}\over{w(z)}}{{p_{n-1}(z)}\over{h_{n-1}z^{n+1}}}\int_S
t^np_n(t)w(t)dt+O\Bigl( {{1}\over{z^2}}\Bigr)\cr
&={{n}\over{z}}+O\Bigl( {{1}\over{z^2}}\Bigr)\qquad (z\rightarrow\infty
).&(4.11)\cr}$$
\noindent By comparing the coefficients of $1/z$ in (4.7) with (4.9),
(4.10) and (4.11), 
we obtain
$$\sum_{k=1}^{4N-2}\alpha_k(n)=\left[\matrix{ n&0\cr
0&1-n-2N\cr}\right],\eqno(4.12)$$
\noindent which leads to (4.4).\par
\rightline{$\square$}\par
\vskip.05in
\indent Let
$$\Phi_n(z)=\left[\matrix{ \sqrt{2\pi i}p_n(z)& -{{i\pi
w(z)p_n(z)+q_n(z)}\over{w(z)\sqrt{2\pi i}}}\cr
{{\sqrt{2\pi i} p_{n-1}(z)}\over{h_{n-1}}}& -{{i\pi
w(z)p_{n-1}(z)+q_{n-1}(z)}\over{w(z)h_{n-1}\sqrt{2\pi
i}}}\cr}\right],\eqno(4.13)$$
\noindent which is a matrix function with entries in ${\bf C}(z)[w]$;
note that $\Phi_n$ also depends upon the $\delta_j$.\par

\vskip.05in

\noindent {\bf Lemma 4.2} {\sl The functions $\Phi_n$ satisfy\par
\indent (i) the basic differential equation
$${{d\Phi_n(z)}\over{dz}}=A_n(z)\Phi_n(z),\eqno(4.14)$$
\indent (ii) the deformation equation
$${{\partial \Phi_n}\over{\partial
\delta_j}}=-{{\alpha_j(n)}\over{z-\delta_j}}\Phi_n(z),\eqno(4.15)$$
\indent (iii) and the recurrence relation 
$\Phi_{n+1}(z)=V_n(z)\Phi_n(z);$\par
\indent (iv) moreover, $\Phi_n$ is invertible since 
$\det \Phi_n(z)=1/w(z).$}\par 
\vskip.05in
\noindent {\bf Proof.} (i) We can write
$$\Phi_n(z)=Y_n(z)\left[\matrix{\sqrt{2\pi i}&0\cr 0&1/(w(z)\sqrt{2\pi
i})\cr}\right],\eqno(4.16)$$
\noindent and then the property (i) follows from (4.1).\par
\indent (ii) This follows
from (i) by standard results in the theory of Fuchsian differential
equations as in [14, 16].\par
\indent (iii) The recurrence relation from Proposition 3.2(i).\par
\indent (iv) Given (iii), this identity follows from 
Proposition 3.2(ii).\par
\rightline{$\square$}\par
\vskip.05in
\noindent Lemma 4.2 states several properties that the $\Phi_n$ satisfy
simultaneously, and hence generates several consistency conditions. By
taking (i), (ii)  and (iii) pairwise, we obtain three Lax pairs, which
we state in the following three propositions.

\vskip.05in
\noindent {\bf Proposition 4.3} {\sl The residue matrices satisfy
Schlesinger's equations
$${{\partial
\alpha_k(n)}\over{\partial\delta_j}}={{[\alpha_j(n),\alpha_k(n)]}\over{
\delta_j-\delta_k}}\qquad (j\neq k)\eqno(4.17)$$
\noindent and}
$${{\partial
\alpha_j(n)}\over{\partial\delta_j}}=-\sum_{k=1;
j\neq k}^{4N-2}{{[\alpha_j(n),\alpha_k(n)]}\over{
\delta_j-\delta_k}}.\eqno(4.18)$$
\par
\vskip.05in
\noindent {\bf Proof.} We can express the consistency condition
${{\partial^2\Phi_n(z)}\over {\partial \delta_j\partial z}}=
{{\partial^2\Phi_n(z)}\over {\partial z\partial \delta_j}}$ as the Lax
pair
$${{\partial A_n(z)}\over {\partial
\delta_j}}-A_n(z){{\alpha_j(n)}\over{z-\delta_j}}={{\alpha_j(n)}\over{
(z-\delta_j)^2}}-{{\alpha_j(n)A_n(z)}\over{z-\delta_j}}\eqno(4.19)$$
\noindent and then one can simplify the resulting system of
differential equations. See [14, 19].\par
\rightline{$\square$}\par
\vskip.05in
\noindent {\bf Proposition 4.4} {\sl The basic differential equation
(4.14)
and the recurrence relation in Lemma 4.2 are consistent, so}
$$A_{n+1}(z)V_n(z)-V_n(z)A_n(z)=\left[\matrix{1&0\cr
0&0\cr}\right].\eqno(4.20)$$
\vskip.05in
\noindent {\bf Proof.} The Lax pair associated with the
these conditions gives
$$A_{n+1}(z)\Phi_{n+1}(z)={{d}\over{dz}}\Phi_{n+1}(z)={{d}\over{dz}}
\Bigl( V_n(z)\Phi_n(z)\Bigr).\eqno(4.21)$$
\rightline{$\square$}\par
\vskip.05in
\noindent {\bf Proposition 4.5} {\sl (i) The deformation equation
(4.15) and
the recurrence relation in Lemma 4.2 are consistent, so 
$$-{{\alpha_j(n+1)}\over{z-\delta_j}}V_n(z)+V_n(z)
{{\alpha_j(n)}\over{z-\delta_j}}={{\partial V_n(z)}\over{\partial
\delta_j}}.\eqno(4.22)$$
\indent (ii) In particular, the $(1,2)$ entry satisfies}
$${{\partial }\over{\partial\delta_j}}\log
h_n=-h_n^{-1}\alpha_j(n)_{12}.\eqno(4.23)$$
\rightline{$\square$}\par
\vskip.05in
\noindent {\bf Proof.} (i) This is the Lax pair associated with
Lemma 4.2.\par
\indent (ii) By letting $z\rightarrow\infty $ in (4.22), we deduce
$$-\alpha_j(n+1)\left[\matrix{ 1&0\cr 0&0\cr}\right] +
\left[\matrix{ 1&0\cr 0&0\cr}\right] \alpha_j(n)
=\left[\matrix{-{{\partial b_{n+1}}\over{\partial \delta_j}}&
-{{\partial h_n}\over{\partial\delta_j}}\cr
-{{1}\over{h_n^2}}{{\partial h_n}\over{\partial\delta_j}}&
0\cr}\right]\eqno(4.24)$$
\noindent which implies that $\alpha_j(n)_{12}=-{{\partial
h_n}\over{\partial\delta_j}}.$\par
\rightline{$\square$}\par
\vskip.05in
\noindent {\bf 5. The tau function}\par
\vskip.05in
\indent We
introduce the differential $1$-form on ${\bf C}^{4N-2}\setminus \{
{\hbox{diagonals}}\}$ by 
$$\Omega_n =\sum_{j,k=1; j\neq k}^{4N-2}
{\hbox{trace}}\Bigl(
{{\alpha_j(n)\alpha_k(n)}\over{\delta_j-\delta_k}}\Bigr)
d\delta_j.\eqno(5.1)$$
\vskip.05in
\noindent {\bf Theorem 5.1} {\sl The Hankel determinant $D_n$ gives the tau
function, so }
$$\Omega_n =d\log D_n.\eqno(5.2)$$
 \vskip.05in
\noindent {\bf Proof.} By Proposition 4.3 and results of Jimbo {\sl et
al} [20], $\Omega_n$ is an exact differential form, so $d\Omega_n=0$; hence
there exists a function $\tau_n$ such that $d\log \tau_n=\Omega_n$,
and so we proceed to identify $\tau_n$. By Lemma 3.1(iii), we have 
$\log h_n=\log D_{n+1}/D_n$, hence we consider
$$\Omega_{n+1}-\Omega_n=\sum_{j\neq k: j,k=1}^{4N-2}
{\hbox{trace}}\Bigl(
{{\alpha_{j}(n+1)\alpha_k(n+1)-
\alpha_j(n)\alpha_k(n)}\over{\delta_j-\delta_k}}\Bigr)d\delta_j
\eqno(5.3)$$
\noindent where by Proposition 4.5(i)
$\alpha_j(n+1)=V_n(\delta_j)\alpha_j(n)V_n(\delta_j)^{-1}$ so
$$\eqalignno{{\hbox{trace}}\bigl(&
\alpha_j(n+1)\alpha_k(n+1)-\alpha_j(n)\alpha_k(n)\bigr)
\cr
&={\hbox{trace}}\bigl(
\alpha_j(n)V_n(\delta_j)^{-1}V_n(\delta_k)\alpha_k(n)V_n
(\delta_k)^{-1}V_n(\delta_j)-\alpha_j(n)\alpha_k(n)\bigr).&(5.4)\cr}$$
\indent We have
$$V_n(\delta_j)^{-1}V_n(\delta_k)=\left[\matrix{ 1&0\cr
{{\delta_j-\delta_k}\over{h_n}}& 1\cr}\right]\eqno(5.5)$$
\noindent so by direct calculation
$$\eqalignno{\Omega_{n+1}-\Omega_n=
\sum_{j\neq k: j,k=1}^{4N-2} &\Bigl\{h_n^{-1}\alpha_j(n)_{12}\bigl(
\alpha_k(n)_{11}-\alpha_k(n)_{22}\bigr)\cr
& +h_n^{-1}\alpha_k(n)_{12}\bigl(
\alpha_j(n)_{22}-\alpha_j(n)_{11}\bigr)\cr
&-h_n^{-2}(\delta_j-\delta_k)\alpha_j(n)_{12}\alpha_k(n)_{12}\Bigr\}
d\delta_j&(5.6)\cr}$$
In this sum we have taken $j\neq k$, but the expression is unchanged if
we include the corresponding terms for $j=k$; hence the coefficient of
$d\delta_j$ is
$$\eqalignno{\alpha_j(n)_{12}&\sum_{k=1}^{4N-2}h_n^{-1}(
\alpha_k(n)_{11}-\alpha_k(n)_{22})
-(\alpha_j(n)_{11}-\alpha_j(n)_{22})\sum_{k=1}^{4N-2}h_n\alpha_k(n)_{12}\cr
&\quad
-{{\delta_j\alpha_j(n)_{12}}\over{h_n^2}}\sum_{k=1}^{4N-2}
\alpha_k(n)_{12}+{{\alpha_j(n)_{12}}\over{h_n^2}}
\sum_{k=1}^{4N-2}\delta_k\alpha_k(n)_{12}.
&(5.7)\cr}$$
\noindent We use Lemma 4.1 to reduce this to
$-h_n^{-1}\alpha_j(n)_{12}$, so 
$$\eqalignno{\Omega_{n+1}-\Omega_n&
=-\sum_{j=1}^{4N-2} h_n^{-1}\alpha_j(n)_{12}d\delta_j \cr
&=\sum_{j=1}^{4N-2} {{\partial }\over{\partial\delta_j}}\log h_n
d\delta_j.&(5.8)\cr}$$ 
\noindent Hence $\Omega_n=d\sum_{k=0}^{n-1}\log h_k.$\par
\rightline{$\square$}\par
\vskip.05in

\indent Following [19], we interpret (5.1) in terms of integrable systems and Hamiltonian mechanics. Let $M=M_2({\bf R})^{4N-2}$ be the product space of matrices,
and let $G=GL_2({\bf R})$ act on $M$ by conjugating each matrix in the
list
$$(X_1, \dots , X_n)\mapsto (UX_1U^{-1}, \dots , UX_nU^{-1}).$$
 The Lie algebra ${\hbox{g}}$ of $G$ has dual ${\hbox{g}}^*$, and
for each $\xi\in {\hbox{g}}^*$ the symplectic structure at $\xi$ on 
${\hbox{g}}\times {\hbox{g}}$ is given by $\omega_\xi (X,Y)=\xi
([X,Y]).$ Given 
$$A(z)=\sum_{k=1}^{2N-2}{{\alpha_k}\over{z-\delta_k}}\eqno(5.9)$$
\noindent as in (4.2), we introduce 
$$\omega (X,Y)=\sum_{k=1}^{2N-2} {\hbox{trace}}\bigl(
\alpha_k[X_k,Y_k]\bigr)\eqno(5.10)$$
\noindent for $X=(X_k)_{k=1}^{2N-2}$ and  $Y=(Y_k)_{k=1}^{2N-2}$ in
${\hbox{g}}^{2N-2}.$ Given $f,g:M\rightarrow {\bf C}$, their Poisson
bracket is $\{ f,g\}=X_f(g)$, and the corresponding vector field
satisfies $X_{\{ f,g\}} =[X_f,X_g].$ The spectral curve of $A(z)$ is the algebraic variety
$$\Sigma_A=\Bigl\{ (z,w)\in {\bf C}^2: \det \bigl(
wI-A(z)\bigr)=0\Bigr\}.\eqno(5.11)$$
\indent As suggested by (5.1), we introduce the Hamiltonian
$H_j:M\rightarrow {\bf C}$ by
$$H_j=\sum_{k:k\neq j}{\hbox{trace}}\Bigl( 
{{\alpha_j(n)\alpha_k(n)}\over{\delta_j-\delta_k}}\Bigr)\eqno(5.12)$$
\noindent so that ${{\partial }\over{\partial \delta_j}}\log \tau
(\delta )=H_j$. We observe that $H_j$ is a polynomial in the
entries of $\alpha_j(n)$ and $\alpha_k(n)$, and is a rational function
of $\delta_j$ and $\delta_k$. To lighten the notation, we
temporarily suppress
the variable $n$.\par
\vskip.05in
\noindent {\bf Proposition 5.2} {\sl (i) The Hamiltonian $H_j$ gives a
vector field $(X_{H_j}^{(k)})_{k=1}^{2N-2}$ which is associated with the
differential equation}
$${{dA}\over{dt}}=\Bigl[ A,
{{\alpha_j}\over{z-\delta_j}}\Bigr].\eqno(5.13)$$
\indent {\sl (ii) The Poisson brackets of the flows commute, so that 
$\{ H_j, H_k\}=0.$}\par
\indent {\sl (iii) Under this flow, the spectral curve of $A$ is
invariant.}\par
\vskip.05in
\noindent {\bf Proof.} (i) For each $Y=(Y_k)_{k=1}^{2N-2}$, we introduce a flow on
$M$ by $\dot \alpha_k=[Y_k,\alpha_k]$. We can differentiate $H_j$ in
the direction of $Y$ and obtain
$$Y(H_j)=\sum_{k:k\neq j}{{{\hbox{trace}}\,
([Y_k,\alpha_k]\alpha_j)}\over{\delta_j-\delta_k}}+
{{{\hbox{trace}}\,
(\alpha_k[Y_k,\alpha_j])}\over{\delta_j-\delta_k}}\eqno(5.14)$$
\noindent With $H_k$, we associate the Hamiltonian vector field
$X_{H_j}=(X_{H_j}^{(k)})_{k=1}^{2N-2}$ such that 
$$Y(H_j)=\omega (X_{H_j}, Y)=\sum_{k=1}^{2N-2} {\hbox{trace}}\bigl(
\alpha_k[(X_{H_j}^{(k)},Y_k]\bigr).\eqno(5.15)$$
\noindent We deduce that 
$$\eqalignno{X_{H_j}^{(k)}&={{\alpha_j}\over{\delta_j-\delta_k}}\qquad (k\neq j)
&(5.16)\cr
X_{H_j}^{(j)}&=\sum_{k:k\neq j}
{{\alpha_k}\over{\delta_j-\delta_k}}.&(5.17)\cr}$$
\noindent It is then a simple calculation to check that $\dot
\alpha_k=[X_{H_j}^{(k)},\alpha_k]$ extends to give (5.13) for $(d/dt) A(z)$.\par
\indent (ii) Given the vector fields $(X_{H_k}^{(j)})_{j}$
corresponding to $H_k$ and  $(X_{H_\ell}^{(j)})_{j}$
corresponding to $H_\ell$ from (5.12), one can compute
$$\{ H_k, H_\ell\} =\sum_j{\hbox{trace}}\bigl( [X_{H_k}^{(j)},
A_j]X_{H_\ell}^{(j)}\bigr)\eqno(5.18)$$
\noindent and reduce the expression to zero by an elementary
calculation.\par
\indent (iii) One can check that for each positive integer $m$,
the ${{d}\over{dt}}{\hbox{trace}}A(z)^m=0$, and hence $\det (wI-A(z))$
is invariant under the flow.\par 
\rightline{$\square$}\par
\vskip.05in
\noindent {\bf 6. Orthogonal polynomials on a single interval}\par
\vskip.05in
\noindent In this section we consider the Chebyshev polynomials.\par
\indent $\bullet$ First suppose that $v=0$. Then the corresponding
equilibrium distribution is the Chebyshev distribution $(1/\pi
)(1-x^{2})^{-1/2}$. In this case, the orthogonal polynomials are the Chebyshev
polynomials of the first kind and, unsurprisingly, our results reduce
to those of Chen and Its [8].\par 
\vskip.05in
\indent $\bullet$ For $a<b$, let 
$$v(z)={{8}\over{(b-a)^2}}\Bigl( z-{{a+b}\over{2}}\Bigr)^2,\eqno(6.1)$$
\noindent so, by standard results used in random matrix theory [26], 
the equilibrium measure is the semicircular law on $[a,b]$, as given by 
$$\rho (dx)={{8}\over{\pi (b-a)^2}}\sqrt{ (b-x)(x-a)}\,{\bf
I}_{[a,b]}(x)dx\eqno(6.2)$$
\vskip.05in

\noindent {\bf Proposition 6.1} {\sl The tau function for the semicircular
distribution is} 
$$\tau (a,b)=(a-b)^{(2n^2+2n+1)/4}e^{n(n+2)(a-b)^2/32}.\eqno(6.3)$$
\vskip.05in
\noindent {\bf Proof.}  Let $U_n$ be the Chebyshev polynomial of the second kind of
degree $n$, which satisfies
$$U_n(\cos\theta )={{\sin (n+1)\theta}\over{\sin\theta}},\eqno(6.4)$$
\noindent and let 
$$p_n(x)=2^{-2n}(b-a)^nU_n\Bigl(
{{2x-(a+b)}\over{b-a}}\Bigr)\eqno(6.5)$$
\noindent which is monic and of degree $n$, and the $p_n$ are
orthogonal with respect to the measure $\rho$.
\noindent By elementary calculations involving trigonometric functions,
one can show that $h_n=2^{-4n}(b-a)^{2n}$ and 
$$A_n(x)={{1}\over{(x-b)(x-a)}}\left[\matrix{ n(x-(a+b)/2)&
-(n+1)(b-a)^2h_{n-1}/8\cr
n(b-a)^2/(2h_{n-1})& -(n+1)(x-(a+b)/2)\cr}\right],\eqno(6.6)$$
\noindent which has poles at $a$ and $b$, as expected. One
verifies that 
$$\Omega_n=\Bigl( {{n^2+(n+1)^2}\over{4}}\Bigr)
{{da-db}\over{a-b}}+{{n(n+2)}\over{16}}(a-b)(da-db)\eqno(6.7)$$
\noindent so that (6.3) follows by integration.
\rightline{$\square$}\par
\vskip.05in
\noindent {\bf 7. Painlev\'e equations for pairs of intervals}\par
\vskip.05in
\noindent Akhiezer considered a generalization of the Chebyshev
polynomials to the pair of intervals $[-1, \alpha ]\cup [\beta , 1]$,
and investigated their properties by conformal mapping. Chen and
Lawrence [9] used the theory of elliptic functions to investigate these
polynomials and in (8.18) expressed the Hankel determinant in terms of
Jacobi's elliptic theta functions. In this section we obtain the 
differential equation where $S$ is two intervals, and obtain a differential equation for the
endpoints that is related to the one from [9].\par
\indent Let $v$ be a polynomial of degree $2N\geq 4$ such that
$S=[\delta_1, \delta_2]\cup [\delta_3, \delta_4]$. There exists a
M\"obius transformation $\varphi $ such that $\varphi (\delta_1)=0$,
$\varphi (\delta_2)=1$ and $\varphi (\delta_4)=\infty$; then we let
$t=\varphi (\delta_3).$ Having fixed three of the endpoints, we can
introduce the differential equations from section 4 that describe the
effect of varying the endpoint $t$, namely
$${{d}\over{dx}}\Phi
(x)=\Bigl({{\alpha_0}\over{x}}+{{\alpha_1}\over{x-1}}+{{\alpha_t}\over{
x-t}}\Bigr) \Phi\eqno(7.1)$$
\noindent and
$${{\partial \Phi}\over{\partial t}}={{-\alpha_t}\over{x-t}}\Phi
.\eqno(7.2)$$
\noindent Let $A(x,t)$ be the matrix
$(\alpha_0/x+\alpha_1/(x-1)+\alpha_t/(x-t))$ and let $A(x,t)_{12}$ be
its top right entry. 
Then we introduce $x=\lambda (t)$ such that
$A(x,t)_{12}=0;$ then by [18, p. 1333], the corresponding Schlesinger equations give a version
of the nonlinear Painlev\'e equation $P_{VI}$ in terms of
$\lambda$, namely
$$\eqalignno{{{d^2\lambda }\over{dt^2}}+&\Bigl({{1}\over{t}}+
{{1}\over{t-1}}+{{1}
\over{\lambda -t}}\Bigr) {{d\lambda}\over{dt}}-{{1}\over{2}}
\Bigl({{1}\over{\lambda}}+{{1}\over{\lambda
-1}}+{{1}\over{\lambda -t}}\Bigr) \Bigl( {{d\lambda}\over{dt}}\Bigr)^2
\cr
&={{1}\over{2}}{{\lambda (\lambda -1)(\lambda
-t)}\over{t^2(t-1)^2}}\Bigl( k_\infty
-{{k_0t}\over{\lambda^2}}+{{k_1(t-1)}\over{(\lambda -1)^2}}-
{{(k_t-1)t(t-1)}\over{(\lambda -t)^2}}\Bigr).&(7.3)\cr}$$
\indent The Hamiltonian and tau function
satisfy
$$H_t= {{d}\over{dt}}\log\tau ={\hbox{trace}}\Bigl(
{{\alpha_0\alpha_t}\over{t}}+{{\alpha_1\alpha_t}\over{t-1}}\Bigr).\eqno
(7.4)$$
\vskip.05in

\noindent Having transformed $S$ to $[0,1]\cup [t,\infty ]$ we can lift
this to the portions of the real axis that are covered by the elliptic 
curve ${\cal E}=\{ (\lambda ,w):w^2=4\lambda (\lambda -1)(\lambda -t)\}$ which has parameters 
$\lambda ={\cal P}(u/2)$ and $w={\cal P}'(u/2)$ in terms of
 Weierstrass's function ${\cal P}$ with
$e_1=t$, $e_2=1$ and $e_3=0$. Hence we transform to the 
dependent variables
$$u=\int_0^\lambda {{ds}\over{\sqrt{s(s-1)(s-t)}}},\eqno(7.5)$$
\noindent  and make the substitution 
$u=u(\lambda (t),t)$.  Fuchs [16] observed that 
$$\eqalignno{{{d^2u}\over{dt^2}}+&{{2t-1}\over{t(t-1)}}
{{du}\over{dt}}+{{u}\over{4t(t
-1)}}\cr
&={{\sqrt{\lambda (\lambda -1)(\lambda -t)}}\over{2t^2(t-1)^2}}\Bigl[
k_\infty -{{k_0t}\over{\lambda^2}}+{{k_1(t-1)}\over{(\lambda -1)^2}}-
{{k_tt(t-1)}\over {(\lambda -t)^2}}\Bigr].&(7.6)\cr}$$
\noindent To solve this in a special case, we introduce the complete elliptic integral
$$K(m)=\int_0^{\pi /2} {{d\theta}\over{\sqrt{1-m^2\sin^2\theta}}}.$$ 
\noindent By comparing terms on power series, one can recover the
following result.\par
\vskip.05in
\noindent {\bf Proposition 7.1} (Poincar\'e) {\sl Suppose that
$k_0=k_1=k_t=k_\infty =0$ and that $u(t)=c_1K(\sqrt {t})+c_2K(\sqrt
{1-t})$ for constants $c_1$ and $c_2$ . Then $u$ satisfies Legendre's
equation
$$t(t-1){{d^2u}\over{dt^2}}+(2t-1){{du}\over{dt}}+{{u}\over{4}}=0,\eqno(7.7)$$
\noindent so $\lambda (t)={\cal P}(u(t)/2; 0,1,t)$ gives a solution of 
$P_{VI}$.}\par
\vskip.05in

\noindent {\bf 8. Kernels associated with Schlesinger's equations}\par 
\vskip.05in
\noindent In this section, we introduce kernels that are associated
with Schlesinger's equations, and then factorize them in terms of Hankel
operators. First we let $\nu_j=-2^{-1}{\hbox{trace}}\,\alpha_j(n)$ and
observe that $\nu_j$ does not depend upon $n$. Indeed, by
multiplying (4.22) by $V_n^{-1}$, one deduces that
${\hbox{trace}}\,A_{n+1}(z)={\hbox{trace}}\,A_n(z)$, and since
${\hbox{trace}}\,\alpha_j(n)=\lim_{z\rightarrow \delta_j}
(z-\delta_j){\hbox{trace}}\,A_n(z),$ we deduce that
${\hbox{trace}}\, \alpha_j(n)$ is constant with respect to $n$. By
(4.12), we have $\sum_{j=1}^{4N-2} {\hbox{trace}}\, \alpha_j(n)=1-2N.$
Now, given $\Phi_n$ as in (4.13), let
$$\Psi_n(z)=\prod_{j=1}^{4N-2}
(z-\delta_j)^{\nu_j}\Phi_n(z).\eqno(8.1)$$
\noindent We next introduce the matrix valued kernel
$$M_n(z, \zeta )={{\Psi_n(z)^\dagger J\Psi_n(\zeta )}\over
{-2\pi i(z-\zeta )}};\eqno(8.2)$$
\noindent we aim to show that $M_n$ is positive definite as an integral
operator on $L^2(S)$, and we observe that this property does not 
change if we introduce weights on $S$.\par 
\vskip.05in
\noindent{\bf Proposition 8.1} {\sl Let $E_{n}(z,\zeta )$ be the kernel of
the orthogonal projection onto ${\hbox{span}}\{ x^j:j=0, \dots, n-1\}$
in $L^2(\rho ).$ Then the top left entry of $M_n(z, \zeta )$ equals}
$$M_n(z, \zeta )_{11} =
{{h_{n}}\over{h_{n-1}}}\prod_{j=1}^{4N-2}(z-\delta_j)^{\nu_j}\prod_{j=
1}^{4N-2}(\zeta-\delta_j)^{\nu_j}E_{n}(z, \zeta ).\eqno(8.3)$$

\vskip.05in
\noindent {\bf Proof.} The Christoffel--Darboux formula gives
$$E_n(z, \zeta )={{p_n(z)p_{n-1}(\zeta )-p_{n-1}(z)p_n(\zeta )}\over{
h_n (z-\zeta )}}.\eqno(8.4)$$
\noindent One can find $\Psi_n(z)^\dagger J\Psi_n(\zeta )$ by direct
calculation, and compare with this.\par
\rightline{$\square$}\par
\vskip.05in
\indent Let $\beta_j(n)=\alpha_j(n)+\nu_jI_2$, which has zero 
trace. Furthermore, if $\Phi_n$ is a
solution of the basic differential equation (4.14), then
$${{d}\over{dz}}\Psi_n(z)=B_n(z)\Psi_n(z)\eqno(8.5)$$
\noindent where
$$B_n(z)=\sum_{j=1}^{4N-2} {{\beta_j(n)}\over{z-\delta_j}}.\eqno(8.6)$$
\noindent We pause to note an existence result for solutions of the
matrix system (8.5).\par
\vskip.05in
\noindent {\bf Lemma 8.2} {\sl Suppose that $\beta_j(n)$ has eigenvalues $\pm
\kappa_j(n)$ where $2\kappa_j(n)$ is not an integer. Then on a
neighbourhood of $\delta_j$, there
exists an analytic matrix function $\Xi_{n,j}$ such that
$$\Psi_n(z)=\Xi_{n,j}(z)(z-\delta_j)^{\beta_j(n)}\eqno(8.7)$$
\noindent satisfies (8.5).}\par

\vskip.05in
\noindent {\bf Proof.} This follows from Turrittin's theorem; see
[3].\par
\rightline{$\square$}\par
\vskip.05in

\indent For notational simplicity, we consider the interval $(\delta_1,
\delta_2)$ and assume that $\delta_1=0$ and $1<\delta_2$; the general
case follows by scaling. For a continuous function $\phi :(0,
1)\rightarrow {\bf R}^{8N-6}$, the Hankel operator $\Gamma_\phi : 
L^2((0, 1); dy/y; {\bf R})\rightarrow L^2((0, 1); dy/y; {\bf
R}^{8N-6})$ is given by
$$\Gamma_\phi f(x)=\int_0^1 \phi (xy)f(y){{dy}\over{y}}.\eqno(8.8)$$
\vskip.05in
\indent Since $\beta_k(n)$ has zero trace, the matrix
$(\delta_1-\delta_k)J\beta_k(n)$ is real symmetric and hence is
congruent to either
$$\sigma_k=\pm\left[\matrix{1&0\cr 0&1\cr}\right], \quad
\pm\left[\matrix{1&0\cr 0&0\cr}\right],\quad \left[\matrix{1&0\cr
0&-1\cr}\right],\quad \left[\matrix{0&0\cr 0&0\cr}\right];\eqno(8.9)$$
\noindent let $\sigma ={\hbox{diagonal}}[\sigma_k]_{k=2}^{4N-2}$ be the
block diagonal sum of these matrices.  
\vskip.05in

\noindent {\bf Theorem 8.3} {\sl (i) Let $\beta_1(n)$ be as in
Lemma 8.2. Then there exists $Z_n$, a
$2\times 1$ real vector solution of
(8.5) such that $Z_n(x)\rightarrow 0$ as $x\rightarrow
\delta_1$.\par
\indent (ii) The integral operator on $L^2((0,1); dx/x)$ with kernel 
$$K_n(z, \zeta )={{\sqrt{z\zeta} Z_n(\zeta )^\dagger
JZ_n(z)}\over{z-\zeta}}\eqno(8.10)$$
\noindent is of trace class;  moreover, 
there exists a real vector Hankel operator $\Gamma_{\psi_n}$ 
on\par
\noindent $L^2((0, 1), dy/y)$ such that 
$$K_n=\Gamma_{\psi_n}^\dagger \sigma \Gamma_{\psi_n}.\eqno(8.11)$$
\indent (iii) If $\sigma\geq 0$, then  $K_n\geq 0$.}\par 
\vskip.05in
\noindent {\bf Proof.} (i) There exists an invertible constant $2\times 2$
matrix $S_n$ such that
$$S_nz^{\beta_1(n)}S_n^{-1}=\left[\matrix{z^{\kappa_1(n)}&0\cr
0&z^{-\kappa_1(n)}\cr}\right].\eqno(8.12)$$
\noindent where $\kappa_1(n)>0$. Hence by Lemma 8.2, there exists a constant
$2\times 1$ matrix $C$
such that $Z_n(z)=\Psi_n(z)C$ is a solution of (8.5), and
$Z_n(z)= O(\vert z-\delta_1\vert^{\kappa_1(n)})$ as 
$z\rightarrow \delta_1$.\par
\indent (ii) Hence we can introduce
$K_n$ by (8.10), and next we prove that the kernel satisfies  
$$\Bigl(x{{\partial}\over{\partial
x}}+y{{\partial}\over{\partial y}}\Bigr)
K_n(x,y)=\sum_{k=2}^{4N-2}
{{-\delta_k\sqrt{xy}}\over{(x-\delta_k)(y-\delta_k)}}Z_n(y)^\dagger J\beta_k
(n)Z_n(x).\eqno(8.13)$$ 
\noindent First note that by homogeneity $(x{{\partial }\over{\partial
x}}+y{{\partial}\over{\partial y}})\sqrt{xy}/(x-y)=0$. Since the $\beta_k(n)$ have zero trace, we have
$J\beta_k(n)+\beta_k(n)^\dagger J=0$ and hence the differential equation gives
$$\eqalignno{\Bigl(x{{\partial}\over{\partial
x}}&+y{{\partial}\over{\partial
y}}\Bigr)Z_n(y)^\dagger JZ_n(x)\cr
&=Z_n(y)^\dagger B_n(y)^\dagger JZ_n(x)+Z_n(y)^\dagger 
JB_n(x)Z_n(x)\cr
&=\sum_{k=2}^{4N-2}
\sqrt{xy}Z_n(y)^\dagger J\beta_k(n)Z_n(x)\Bigl(
{{x}\over{x-\delta_k}}-{{y}\over{y-\delta_k}}\Bigr);
&(8.14)\cr}$$ 
\noindent on dividing by $x-y$, we obtain
$$\Bigl(x{{\partial}\over{\partial
x}}+y{{\partial}\over{\partial
y}}\Bigr){{\sqrt{xy}Z_n(y)^\dagger JZ_n(x)}\over{x-y}}=
\sum_{k=2}^{4N-2}
{{\delta_k\sqrt{xy}}\over{(x-\delta_k)(y-\delta_k)}}Z_n(y)^\dagger J\beta_k
(n)Z_n(x)\eqno(8.15)$$ 
\noindent as in (8.13).\par
\indent Noting the shape of the final 
factor in (8.12), we choose 
$$\phi_n
(x)={\hbox{column}}\Bigl[{{\sqrt{x}Z_n(x)}\over{x-\delta_k}}\Bigr]_{k=2,
\dots , 4N-2}\eqno(8.16)$$
\noindent which has a $2\times 1$ entry for each endpoint $\delta_k$ of $S$
after $\delta_1$, and the block diagonal matrix
$$\beta (n)={\hbox{diagonal}}\Bigl[ -\delta_kJ\beta_k(n)\Bigr]_{k=2,
\dots , 4N-2}\eqno(8.17)$$
\noindent with $2\times 2$ blocks, and we consider
$$\tilde K_n(x,y)=\int_0^1 \phi_n (yz)^\dagger\beta (n)\phi_n
(zx){{dz}\over{z}}.\eqno(8.18)$$
\noindent First note that since $\kappa_1(n)>0$, we have $\tilde
K(x,y)\rightarrow 0$ as $x,y\rightarrow 0$. Then 
$$\eqalignno{ \Bigl( x{{\partial}\over{\partial
x}}+y{{\partial}\over{\partial y}}\Bigr)\tilde  K_n(x,y)&=\int_0^1 \bigl(
y\phi_n'(yz)^\dagger\beta (n)\phi_n (zx)+x\phi_n (yz)^\dagger\beta
(n)\phi_n'(xy)\bigr)\,
dz\cr
&=\phi_n (y)^\dagger\beta (n)\phi_n (x)-\phi_n (0)^\dagger\beta (n)\phi_n
(0).&(8.19)\cr}$$ 
\noindent We have $\phi_n (0)=0$, so
$$K_n(x,y)=\tilde K_n(x,y)+\xi (x/y)\eqno(8.20)$$
\noindent for some function $\xi$. But
$Z_n(z)/(z-\delta_1)^{\kappa_1(n)}$ is analytic on a neighbourhood
of $\delta_1$, so it is clear that $K_n(x,y)\rightarrow 0$ and 
$\tilde K_n(x,y)\rightarrow 0$ as $x\rightarrow 0$ or
$y\rightarrow 0$; hence $\xi =0$. 

\indent By the choice of $\sigma$, there exists a block diagonal matrix 
$\gamma (n)$ such that $\gamma
(n)^\dagger\sigma \gamma (n)=\beta (n)$, so we can introduce
$\psi_n(x)=\gamma (n)\phi_n(x)$ such that $\phi_n(x)^\dagger\beta
(n)\phi_n(y)=\psi_n(x)^\dagger \sigma \psi_n(y)$. For this symbol function
$\psi_n$ we have
$$K_n(x,y)=\int_0^1\psi_n(yz)^\dagger
\sigma \psi_n(zx){{dz}\over{z}},\eqno(8.21)$$
\noindent or in terms of Hankel operators
$K_n=\Gamma_{\psi_n}^\dagger\sigma \Gamma_{\psi_n}$. We have 
$$\int_0^1\log{{1}\over{u}}\,\Vert\psi_n (u)\Vert^2\,
{{du}\over{u}}<\infty ,$$
\noindent so $\Gamma_\psi$ is Hilbert--Schmidt and hence $K_n$ is trace
class.\par
\indent (iii) If $\sigma_k\geq 0$ for all $k$ or equivalently 
$\sigma \geq 0$, then $K_n\geq 0$.\par
\rightline{$\square$}\par
\vskip.05in
\noindent {\bf Corollary 8.4} {\sl Suppose that $Z$ is a
$2\times 1$ solution of
$${{d}\over{dx}}Z
=\Bigl({{\beta_0}\over{x}}+{{\beta_1}\over{x-1}}+{{\beta_t}\over{x-t}}
\Bigr)Z\eqno(8.22)$$
\noindent such that the entries satisfy $Z(\bar x)=\bar Z(x)$ and where\par
\indent (i) $\beta_0$ 
is as in Lemma 8.2, and $Z
(x)\rightarrow 0$ as $x\rightarrow 0$;\par
\indent (ii) $J\beta_1$ is positive definite;\par
\indent (iii) $J\beta_t\geq 0$.\par
\noindent Then there exist an invertible real matrix $S$ and a real
diagonal matrix $D$ such that 
$$\psi (x)=\left[\matrix{{{\sqrt{x} SZ(x)}\over{x-1}}\cr {{\sqrt{x}DSZ
(x)}\over{x-t}}\cr}\right]\eqno(8.23)$$
\noindent satisfies $\psi (\bar x)=\overline{\psi (x)}$ and }
$${{\sqrt{xy}Z(y)^\dagger JZ(x)}\over{y-x}}=\int_0^1 \psi (xz)^\dagger \psi
(zy){{dz}\over {z}}.\eqno(8.24)$$
\vskip.01in
\noindent {\bf Proof.} We simultaneously reduce the quadratic
forms associated with $J\beta_1$ and $J\beta_t$, and introduce an
invertible real matrix $S$ such that $J\beta_1=SS^\dagger$ and
$J\beta_t=SD^2S^\dagger$, where $D$ is a real diagonal matrix such that
the diagonal entries $\kappa$ of $D^2$ satisfy $\det
(\beta_1-\kappa\beta_t)=0$. Then we can write
$$\psi (y)^\dagger\psi (x)=\sqrt{xy} Z(y)^\dagger \Bigl(
{{J\beta_1}\over{(x-1)(y-1)}}+{{J\beta_tt}\over{(x-t)(y-t)}}\Bigr)Z
(x).\eqno(8.25)$$
\noindent Now we can follow the proof of Theorem 8.3, and deduce that
$$-\psi (y)^\dagger\psi (x)=\Bigl( x{{\partial}\over{\partial
x}}+y{{\partial }\over{\partial y}}\Bigr) 
{{\sqrt{xy}Z(y)^\dagger JZ(x)}\over{x-y}};\eqno(8.26)$$
\noindent hence we can obtain the result by integrating and using
(i).\par
\rightline{$\square$}\par
\vskip.05in

\noindent {\bf Example 8.5} For the semicircle law on $[a,b]$, as in
Proposition 6.1, we have
$$J\beta_a(n)=\left[ \matrix{{{n(b-a)}\over{2h_{n-1}}}& 
{{2n+1}\over{4}}\cr
 {{2n+1}\over{4}}& {{(n+1)(b-a)h_{n-1}}\over{8}}\cr}\right]\eqno(8.27)$$
\noindent and
$$J\beta_b(n)=\left[ \matrix{ {{-n(b-a)}\over{2h_{n-1}}}& {{2n+1}\over{4}}\cr
{{2n+1}\over{4}}&
{{-(n+1)(b-a)h_{n-1}}\over{8}}\cr}\right],\eqno(8.28)$$
\noindent so that 
$$\det J\beta_a(n)=\det
J\beta_b(n)={{n(n+1)(b-a)^2}\over{16}}-{{(2n+1)^2}\over{16}}.\eqno(8.29)$$
\noindent In particular, when $a=-1$ and $b=1$, the matrices
$J\beta_{-1}(n)$ and $J\beta_1(n)$ are indefinite.\par

\vskip.05in
\noindent {\bf 9. The tau function realised by a linear system}\par
\vskip.05in
\indent In this section, we express the tau function of $K_n$
from Theorem 8.3 as a Fredholm
determinant, and then obtain this from the solution of an integral
equation of Gelfand--Levitan type. The first step is introduce a
scattering function $\psi$ and
then to realise this by a linear system, so that we can solve 
the Gelfand--Levitan equation.\par
\indent  The
differential equation
$${{dZ_n}\over{dx}}=B_n(x)Z_n(x)\eqno(9.1)$$
\noindent has a solution from which we constructed a symbol function 
$$\psi_n(x) ={\hbox{column}}\Bigl[
{{\sqrt{x}\gamma
(n)Z(x)}\over{x-\delta_k}}\Bigr]_{k=2}^{4N-2}.\eqno(9.2)$$
\noindent Suppressing $n$ for simplicity, we
change $x\in (0,1)$ to $t\in (0, \infty )$ by letting
$x=\delta_1+e^{-t}$ and in the new variables write
$$\psi (t)=\sum_{\ell =0}^\infty \chi_\ell e^{-(\kappa_1+\ell
+1/2)t}.\eqno(9.3)$$
\noindent where $\sum_{\ell=0}^\infty \Vert\chi_\ell\Vert<\infty $.
Likewise, we write $\tau (t)$ for $\tau (\delta_1+e^{-t}).$ \par 
\vskip.05in
\indent Let 
$\Omega=\{ z: \Re z\geq 0\}$ be the open right half-plane, let
$$\Psi (x)
=\left[\matrix{0&\psi (x)\cr \psi
(\bar x)^\dagger&0\cr}\right]\eqno(9.4)$$
and extend $\Psi$ to an analytic function 
$\Psi :\Omega\rightarrow M_{8N-5}({\bf C})$ such that $\Psi (x)=\Psi (x)^\dagger$ for $x>0$. Let $\Psi_{(s)}=\Psi (x+2s)$ and
$\Psi^*_{(s)}(x)=\Psi (x+2\bar
s)^\dagger$ and let $\sigma$ be a constant matrix; then let
$K_s=\Gamma_{\Psi_{(s)}^*}\sigma\Gamma_{\Psi_{(s)}}$ be a family of
operators on $L^2(0, \infty )$. \par
\vskip.05in
\noindent {\bf Proposition 9.1} {\sl (i) The $\tau$ function
associated with $K=\Gamma_{\Psi^*}\sigma\Gamma_{\Psi}$ is  
$\tau (2s)=\det (I-K_s)$, which gives an analytic function on
$\Omega$.\par
\indent (ii) Let
$q(s)=-2{{d^2}\over{ds^2}}\log\tau (2s)$. Then $q(s)$ is meromorphic
on $\Omega$, and analytic where $\int_0^\infty x\Vert\Psi
(x+s)\Vert^2dx<1$.\par
\indent (iii) If $0\leq K\leq I$, then $\tau (s)$ is non-negative for
$0<s<\infty$,
increasing and converges to one as $s\rightarrow\infty$.}\par
\vskip.05in
\noindent {\bf Proof.} (i)  The kernel of the Hankel operator 
$\Gamma_{\Psi_{(s)}}$ has a nuclear
expansion 
$$\Gamma_{\Psi_{(s)}}\leftrightarrow \sum_{\ell =0}^\infty 
e^{-(\kappa_1+\ell +1/2)(x+y+2s)}\left[\matrix{0&\chi_\ell \cr
\chi_\ell^\dagger &0\cr}\right]\eqno(9.5)$$
\noindent where 
$\sum_{\ell =0}^\infty \Vert\chi_\ell\Vert \int_0^\infty
e^{-2(\kappa_1+\ell +1/2)(x+\Re s) }dx <\infty$, so the Fredholm determinants are well
defined. As in Schwarz's reflection principle,
$s\mapsto \Psi^*_{(s)}$ is analytic, and $\Gamma_{\Psi_{(s)}}$ is
Hilbert--Schmidt, so $K_s$ is an analytic trace-class valued function on
$\Omega$. Using unitary equivalence, one checks that 
$$\det (I-K_s)=\det (I-P_{(2s, \infty )}K)\qquad (s>0).\eqno(9.6)$$
\indent (ii)  Except on the discrete set of zeros of
$\tau (2s)$, the operator $I-K_s$ is invertible and 
$$q(s)=2{{d}\over{ds}}{\hbox{trace}}\Bigl(
(I-K_s)^{-1}{{dK_s}\over{ds}}\Bigr) .\eqno(9.7)$$
\indent (iii) This follows from (9.6).\par
\rightline{$\square$}\par 

\indent Next, we obtain an alternative formula for $q$
by realising $\Psi$ via a linear system. The technique is suggested
by the inverse scattering transform. Let $H_0={\bf C}^{8N-6}$ be the column vectors, $H=\ell^2$ be
Hilbert sequence space, written as infinite columns,   
and introduce an infinite row of column vectors $C\in \ell^2(H_0)$ by $C=(\chi_\ell
/\Vert\chi_\ell\Vert^{1/2})_{\ell =0}^\infty$ and a column $B\in
\ell^2$ by $B=(\Vert\chi_\ell\Vert^{1/2})_{\ell =0}^\infty$ and the
infinite square matrix $A={\hbox{diagonal}}\, [ \ell +\kappa_1
+1/2]_{\ell
=0}^\infty $. While $A$ is real and symmetric, we will write
$A^\dagger$ in some subsequent formulas, so as to emphasize their
symmetry.\par
\indent In the following result we use the $(8N-5)\times (8N-5)$ 
block matrices 
$$W(x,y)=\left[\matrix{U(x,y)&v(x,y)\cr w(x,y)^\dagger
&z(x,y)\cr}\right],\quad  \Psi (x)
=\left[\matrix{0&\psi (x)\cr \psi
(\bar x)^\dagger&0\cr}\right],\eqno(9.8)$$
\noindent so that $\Psi (\bar x)=\Psi(x)^\dagger$ and the matrix Hamiltonian
$$H(x)=\left[\matrix{U(x,x)\sigma &v(x,x)\cr
w(x,x)^\dagger \sigma &z(x,x)\cr}\right]\eqno(9.9)$$
\noindent where $v,w\in H_0$, $U$ operates upon $H_0$ and $z$ is a
scalar. To simplify the statements of results, we use a special
 non-associative product $\ast$, involving $\sigma$, that is defined by
$$\left[\matrix{U&v\cr w^\dagger &z\cr}\right]\ast
\left[\matrix{0&\psi \cr \psi^\dagger &0\cr}\right]=
\left[\matrix{v\psi^\dagger &U\sigma \psi\cr z\psi^\dagger
&w^\dagger \sigma \psi \cr}\right].\eqno(9.10)$$

\vskip.05in  
\noindent {\bf Theorem 9.2} {\sl (i) The symbol $\psi$ is realised by the
linear system $(-A, B,C)$, so 
$$\psi (t)=Ce^{-tA}B.\eqno(9.11)$$
\indent (ii) There exists a solution of the Gelfand--Levitan equation 
$$W(x,y)+\Psi (x+y)+\int_x^\infty W(x,s)\ast \Psi (s+y)\, ds
=0\qquad (0<x<y) \eqno(9.12)$$
\noindent such that the tau function of Proposition 9.1(i) satisfies}
$${{d}\over{dx}}\log\tau (2x)={\hbox{trace}}\,H(x)\qquad
(x>0).\eqno(9.13)$$
\indent {\sl (iii) Suppose moreover that $\int_0^\infty x\Vert \Psi
(x)\Vert^2 dx<1$.  Then}
$$\Bigl({{\partial^2}\over{\partial x^2}} -{{\partial^2}\over{\partial y^2}}\Bigr)
W(x,y)=-2{{dH}\over{dx}}\,W(x,y).\eqno(9.14)$$
\vskip.05in
\vskip.05in
\noindent {\bf Proof.} (i) This identity follows from (9.3). Since
$\kappa_1+\ell +1/2>0$, the semigroup $e^{-tA}={\hbox{diagonal}}\,
[e^{-t(\kappa_1+\ell +1/2)}]_{\ell =0}^\infty$ consists of trace class
operators, and the integrals in the remainder of the proof are
convergent.\par
\indent (ii) We introduce the observability Gramian
$$Q^\sigma _x=\int_x^\infty e^{-sA^\dagger}C^\dagger \sigma Ce^{-sA}\, ds\qquad
(x>0),\eqno(9.15)$$
\noindent modified to take account of $\sigma$, and the usual controllability Gramian
$$L_x=\int_x^\infty e^{-sA}BB^\dagger e^{-sA^\dagger}\, ds,\eqno(9.16)$$
\noindent both of which define trace class operators on $\ell^2$, and
where $L_x\geq 0$. The controllability operator
$\Xi_x:L^2((0, \infty ); H_0)\rightarrow H$ is 
$$\Xi_xf=\int_x^\infty e^{-tA}Bf(s)\, ds\eqno(9.17)$$
\noindent while the observability operator is $\Theta_x:L^2((0, \infty
); H_0)\rightarrow H$ is 
$$\Theta_xf=\int_x^\infty e^{-sA^\dagger}C^\dagger f(s)\,
ds.\eqno(9.18)$$ 
\noindent Finally, we let
$\psi_{(x)}(s)=\psi (s+2x)$, so that $\psi_{(x)}$ is realised by 
$(-A, e^{-xA}B, Ce^{-xA})$. In terms of these operators, we have the basic identities
$$\Gamma_{\psi_{(x)}}=\Theta_x^\dagger \Xi_x, \qquad
\Gamma^\dagger_{\psi_{(x)}}=\Xi_x^\dagger\Theta_x\eqno(9.19)$$
\noindent while
$$L_x=\Xi_x\Xi_x^\dagger \qquad {\hbox{and}}\qquad
Q^\sigma _x=\Theta_x\sigma \Theta_x^\dagger .\eqno(9.20)$$
\noindent Hence we can rearrange the factors in the Fredholm
determinants 
$$\eqalignno{\det (I-\lambda
\Gamma_{\psi_{(x)}}^\dagger \sigma \Gamma_{\psi_{(x)}})&=
\det (I-\lambda \Xi_x^\dagger \Theta_x\sigma \Theta_x^\dagger \Xi_x)\cr
&=\det (I-\lambda \Xi_x\Xi_x^\dagger \Theta_x\sigma \Theta_x^\dagger )\cr
&=\det (I-\lambda L_xQ^\sigma_x).&(9.21)\cr}$$
\noindent  We deduce that 
$$\eqalignno{\log\tau (2x)&=\log\det
(I-\Gamma_{\psi}^\dagger\sigma \Gamma_{\psi}P_{(2x,\infty )})\cr
&=\log\det (I-\sigma \Gamma_{\psi}P_{(2x, \infty )}\Gamma_{\psi}^\dagger )\cr
&=\log\det (I-\sigma \Gamma_{\psi_{(x)}}\Gamma_{\psi_{(x)}}^\dagger )\cr
&=\log\det (I-\Gamma_{\psi_{(x)}}^\dagger\sigma 
\Gamma_{\psi_{(x)}})\cr
&={\hbox{trace}}\, \log (I-L_xQ^\sigma _x),&(9.22)\cr}$$
\noindent and hence
$$\eqalignno{{{d}\over{dx}}\log\tau (2x)&={\hbox{trace}}\Bigl( 
(I-L_xQ^\sigma_x)^{-1}\bigl( e^{-xA}BB^\dagger
e^{-xA^\dagger}Q^\sigma_x+L_xe^{-xA^\dagger}C^\dagger
\sigma Ce^{-xA}\bigr)\Bigr)\cr
&=B^\dagger
e^{-xA^\dagger}Q^\sigma_x(I-L_xQ^\sigma_x)^{-1}e^{-xA}B\cr
&\qquad +{\hbox{trace}}\,
\sigma Ce^{-xA}(I-L_xQ^\sigma _x)^{-1}L_xe^{-xA^\dagger}C^\dagger 
.&(9.23)\cr}$$
\indent The integral equation
$$\eqalignno{&\left[\matrix{ U(x,y)& v(x,y)\cr w(x,y)^\dagger & z(x,y)\cr}\right]
+\left[\matrix{ 0& \psi (x+y)\cr \psi (x+y)^\dagger& 0\cr}\right]\cr
&+\int_x^\infty \left[ \matrix{ U(x,s)&v(x,s)\cr
w(x,s)^\dagger& z(x,s)\cr}\right]\ast \left[\matrix{ 0&\psi
(s+y)\cr \psi (s+y)^\dagger&
0\cr}\right] ds=0&(9.24)\cr}$$
\noindent reduces to the
 identities
$$\eqalignno{ U(x,y)&=-\int_x^\infty v(x,s)\psi (s+y)^\dagger\, ds,\cr
z(x,y)&=-\int_x^\infty w(x,s)^\dagger\sigma\psi (s+y)\,
ds,&(9.25)\cr}$$
\noindent and the pair of integral equations
$$v(x,y)+\psi (x+y)-\int_x^\infty\int_x^\infty v(x,t)\psi (t+s)^\dagger 
\sigma \psi (s+y)\, dsdt=0\eqno(9.26)$$
\noindent and
$$w(x,y)+\psi (x+y)-\int_x^\infty\int_x^\infty \psi (s+y)\psi
(t+s)^\dagger \sigma w(x,t)\, dtds=0.\eqno(9.27)$$
\noindent To solve these integral equations, we let
$$v(x,y)=-Ce^{-xA}(I-L_xQ^\sigma_x)^{-1}e^{-yA}B\eqno(9.28)$$
\noindent and  
$$w(x,y)=-Ce^{-yA}(I-L_xQ^\sigma_x)^{-1}e^{-xA}B;\eqno(9.29)$$
\noindent then by substituting these into (9.25) we obtain the diagonal
blocks of the solution $W$, namely
$$U(x,y)=Ce^{-xA}(I-L_xQ^\sigma_x)^{-1}L_xe^{-yA^\dagger}C^\dagger
\eqno(9.30)$$
\noindent and 
$$z(x,y)=B^\dagger
e^{-yA^\dagger}Q^\sigma_x(I-L_xQ^\sigma_x)^{-1}e^{-xA}B.\eqno(9.31)$$
\noindent Hence we can identify the trace of the solution (9.9) as 
$$\eqalignno{ {\hbox{trace}}\, H(x)&=
{\hbox{trace}}\, \sigma U(x,x)+z(x,x)\cr
&={\hbox{trace}}\,
\sigma Ce^{-xA}(I-L_xQ^\sigma_x)^{-1}e^{-xA^\dagger}C^\dagger \cr
&\qquad +B^\dagger
e^{-xA^\dagger}Q^\sigma_x(I-L_xQ^\sigma_x)^{-1}e^{-xA}B\cr
&={{d}\over{dx}}\log\tau (2x).&(9.32)\cr}$$
\indent (iii) By integrating by parts, we obtain the identity
 $$\Bigl({{\partial^2}\over{\partial x^2}} 
-{{\partial^2}\over{\partial y^2}}\Bigr)
W(x,y)-2{{dH}\over{dx}}\Psi (x+y)
+\int_x^\infty \Bigl({{\partial^2}\over{\partial x^2}} 
-{{\partial^2}\over{\partial s^2}}\Bigr)W(x,s)\ast \Psi (s+y)\,
ds=0\eqno(9.33)$$
\noindent for $0<x<y$. One can easily verify that the product
 $\ast$ and the
standard matrix multiplication satisfy
$(QW)\ast \Psi=Q(W\ast \Psi )$, hence the formula
$$-2{{dH}\over{dx}}W(x,y)-2{{dH}\over{dx}}\Psi (x+y)-\int_x^\infty 
\Bigl(2{{dH}\over{dx}}W(x,s)\Bigr)\ast \Psi (s+y)\,
ds=0\eqno(9.34)$$
\noindent follows from multiplying (9.33) by $-2{{dH}\over{dx}}$, and
this shows that both $-2{{dH}\over{dx}}W(x,y)$ and 
$({{\partial^2}\over{\partial x^2}} 
-{{\partial^2}\over{\partial y^2}})
W(x,y)$ are solutions of the same integral equation. By
uniqueness of solutions, they are equal.\par
\rightline{$\square$}\par
\vskip.05in

\noindent {\bf 10. Integrability of the tau function of a linear
system}\par
\vskip.05in
\indent Let $q(x) =-2{{d}\over{dx}}{\hbox{trace}} H(x)$ and $\tau$ be as
in (9.32). In this
section we describe the properties of $\tau$ in terms of the algebraic theory of differential
equations [30].\par 
\indent Let ${\bf F}$ be a field (of complex functions) with differential 
${\partial}$ that contains the subfield ${\bf C}$ of
constants and adjoin an
element $h$ to form ${\bf F}(h)$, where either:\par
\indent (i) $h=\int g$ for some $g\in {\bf F}$, so $\partial h=g$;\par
\indent (ii) $h=\exp \int g$ for some $g\in {\bf F}$; or\par
\indent (iii) $h$ is algebraic over ${\bf F}$.\par
\vskip.05in
\noindent {\bf Definition.} Let ${\bf F}_j$ $(j=1, \dots ,n)$ be fields 
with differential $\partial$ that contain
the subfield ${\bf C}$ of constants and suppose that 
$${\bf F}_1\subseteq {\bf F}_2\subseteq\dots \subseteq {\bf
F}_n,\eqno(10.1)$$
\noindent where ${\bf F}_j$ arises from ${\bf F}_{j-1}$ by applying some operation
(i), (ii) or (iii). Then ${\bf F}_n$ is a Liouvillian extension of
${\bf F}_1$.\par

\vskip.05in
\noindent {\bf Example.} The tau function (6.3) belongs to some
Liouvillian extension of ${\bf C}(x)$.\par
\vskip.05in
\noindent {\bf Lemma 10.1} {\sl Let $A_n:{\bf C}^n\rightarrow {\bf
C}^n$, $B_n:{\bf C}^m\rightarrow {\bf C}^n$ and $C_n:{\bf
C}^n\rightarrow {\bf C}^k$ be finite matrices, 
such that $MA_n+A_n^\dagger M$ is positive definite for
some positive definite $M$. Let $\psi_n (t)=Ce^{-tA_n}B$, and define 
the corresponding terms as in the proof of Theorem 9.2. Then $\tau_n (x)$ 
belongs to some 
Liouvillian extension field ${\bf F}$ of 
${\bf F}_0={\bf C}(t, e^{-t\kappa_j}, e^{-t\bar\kappa_j}; j=1, \dots ,
n)$, where  $(\kappa_j)_{j=1}^n$ is the spectrum of $A_n$.}\par
\vskip.05in
\noindent {\bf Proof.} By Lyapunov's 
criterion [7], all of the eigenvalues
$\kappa_j$ of $A_n$ satisfy $\Re \kappa_j>0$, hence $\Vert
e^{-tA_n}\Vert$ is of exponential decay as $t\rightarrow\infty$. By
considering the Jordan canonical form of $A_n$, we obtain matrix
polynomials $p_j(t)$ such that $e^{-tA_n}=\sum_{j=1}^n
p_j(t)e^{-\kappa_jt}$. Observe that ${\bf F}_0$
contains all the entries of
$e^{-tA_n}B_nB_n^\dagger e^{-tA_n^\dagger}$ and $e^{-tA_n^\dagger
}C_n^\dagger\sigma C_ne^{-tA_n}$. The operator $L_x$ is an indefinite integral of
$e^{-tA_n}B_nB_n^\dagger e^{-tA_n^\dagger}$ while the operator $Q_x^\sigma$ is an
indefinite integral of $e^{-tA_n^\dagger}C_n^\dagger
\sigma C_ne^{-tA_n}$, hence $L_x$ and $Q_x^\sigma$  have entries in ${\bf F}_0$; moreover, the entries of $(I-L_xQ_x^\sigma)^{-1}$ are quotients of
determinants with elements in ${\bf F}_0$. Hence by (9.32), 
${{d}\over{dx}}\log \tau_n(2x)$ gives an element of ${\bf F}_0$, 
so $\tau_n (x)$ itself is in a Liouvillian
extension ${\bf F}$ of ${\bf F}_0$.\par
\rightline{$\square$}\par
\vskip.05in
\noindent {\bf Theorem 10.2}  {\sl Let $\psi$ be as in Theorem
9.2.\par
\indent (i) There exists a
sequence of finite rank matrices $(A_n)_{n=1}^\infty$, with 
corresponding tau functions $\tau_n$, such that 
${{d}\over{dt}}\log \tau_n (2t)$
belongs to ${\bf C}(e^{-(\kappa_1+1/2)t}, e^{-t})$ and
$\tau_n(2t)\rightarrow \tau (2t)$ as $n\rightarrow\infty$, uniformly on compact subsets of
$\{ t: \Re  t>0\}$.\par
\indent (ii) Suppose further that $\kappa_1$ is rational. Then
there exists a positive integer $N$ such that ${{d}\over{dt}}\log\tau (2t)$ is periodic
with period $2\pi iN$, and $\tau_n(2t)$ is given by elementary
functions as in (10.4) below.}\par  
\vskip.05in
\noindent {\bf Proof.} (i) We introduce the finite rank matrices
$$A_n={\hbox{diagonal}}\, [\kappa_1+1/2, \kappa_1+3/2, \dots , 
\kappa_1+n +1/2,
0,0, \dots ]\eqno(10.2)$$
\noindent so that $\Vert e^{-tA}-e^{-tA_n}\Vert_{c^1}\leq
e^{-(\kappa_1+n+1)\Re t}/(1-e^{-\Re t}).$ Similarly, we cut down $B$ and $C$  and follow through the
proof of Theorem 9.2 to produce the appropriate choice of $W_n(t,t)$ by
the prescription of (9.8). Evidently, the eigenvalues of $e^{-tA_n}$
have the form $e^{-t(\kappa_1+\ell +1/2)}$ where $\ell =0, \dots , n$. 
We observe that $e^{-t(\kappa_1+\ell +1/2)}$
belongs to ${\bf C}(e^{-(\kappa_1+1/2)t}, e^{-t})$ for all $\ell$, 
so $W_n(t,t)$ likewise belongs to 
${\bf C}(e^{-(\kappa_1+1/2)t}, e^{-t})$.\par
\indent (ii) In this case, the set $\{
m\kappa_1+m/2+n\ell; m,n\in {\bf Z}; \ell =0, 1, 2, \dots \}$ is a
finitely generated subgroup of the rationals, and hence has a
smallest positive element $M/N$, where $M,N\in {\bf N}$ with $M<N$.
Then all the terms $N(\kappa_1+\ell
+1/2)$ are positive integers, so $\exp (-(t+2\pi Ni)A)=\exp (-tA)$  for
all $\Re t>0$, hence $\tau'(2t )/\tau (2t)$ is periodic.\par
\indent By Lemma 10.1, there exists a rational function $r_n$ such that 
$${{d}\over{dt}}\log \tau_n(2t)=r_n(e^{-t/N}).\eqno(10.3)$$
\noindent Suppose for simplicity that $r_n(z)/z$ has only simple
poles; then from the partial fractions decomposition, there exist coefficients $\alpha_j,\beta_j$ and $\gamma_j$
and $b_j,c_j$ such that $b_j^2<c_j$, real poles $a_k$ and a polynomial $q_n(z)$ such
that (10.3) integrates to
$$\eqalignno{\log\tau_n (2t)& =q_n(e^{-t/N})+\sum_k \alpha_k\log\vert
e^{-t/N}-a_k\vert +\sum_j\beta_j\log
(e^{-2t/N}+2b_je^{-t/N}+c_j)\cr
&\quad +\sum_j 
{{\gamma_j}\over{\sqrt{c_j-b_j^2}}}
\tan^{-1}{{e^{-t/N}+b_j}\over{\sqrt{c_j-b_j^2}}}.&(10.4)\cr}$$  
\noindent When $r_n(z)/z$ has higher order poles, one likewise obtains
expressions that are similar but more complicated.\par
\rightline{$\square$}\par   

\indent In view of Theorem 10.2(ii) and the results of [8], it is natural to consider how the
properties of $\tau$ relate to integrability of  Schr\"odinger's equation
$-f''+qf=\lambda f$.\par
\vskip.05in
\noindent {\bf Definition.} Let $q$ be meromorphic on ${\bf C}$. As
in [6,17], we say that $q$ is algebro-geometric if
there exists a non-zero $R:{\bf C}^2\rightarrow {\bf C}\cup\{\infty \}$ such
that $x\mapsto R(x;\lambda )$ is meromorphic, $\lambda\mapsto
R(x;\lambda )$ is a polynomial, and 
$$-R'''+4(q-\lambda )R'+2q'R=0.\eqno(10.5)$$
\vskip.05in
\indent Drach [6] observed that Schr\"odinger's
equation is integrable by quadratures for all $\lambda$ only if $q$ is algebro
geometric. Conversely, if $R(x; \lambda )$ is as above, then
$$f(x)=\sqrt{R(x;\lambda )}\exp \Bigl( -\int
{{dt}\over{R(t;\lambda )}}\Bigr)\eqno(10.6)$$ 
\noindent gives a solution to Schr\"odinger's equation. In [6], Brezhnev catalogues several known special
functions which give integrable forms of Schr\"odinger's equation. 
The theory extends to meromorphic potentials on compact Riemann surfaces 
by [13, p. 235] and [23, p. 1122].\par
\indent The 
following result of Gesztesy and Weikard summarizes
various sufficient conditions for a potential to be
algebro-geometric. The initial hypothesis rules out variants of
Bessel's equation.\par
\vskip.05in
\noindent {\bf Theorem 10.3} [17] {\sl Suppose that $-f''+qf=\lambda f$ has a meromorphic fundamental
solution for each $\lambda$ and that either\par
\indent (i) $q$ is rational and bounded at infinity;\par
\indent (ii) $q$ is elliptic, that is, doubly periodic; or\par
\indent (iii) $q$ is periodic, with purely imaginary period, and
 $q$ is bounded on $\{ z\in {\bf C}: \vert \Re z\vert >r\}$ for
some $r>0$.\par
\indent Then $q$ is algebro-geometric.}\par  
\vskip.05in
\indent We proceed to consider the cases (i),(ii) and (iii) of Theorem 10.3, 
the linear systems $(-A,B,C)$ that give
rise to them, and the corresponding $\tau$ functions.\par
\vskip.05in
\noindent {\bf Proposition 10.4} {\sl Suppose that $q$ satisfies
Theorem 10.3(i). Then $f$ has a rational Laplace
transform and hence is the transfer function of a linear system
$(-A_n, B,C)$ with a finite matrix $A_n$.}\par
\vskip.05in
\noindent {\bf Proof.} By a theorem of Halphen [17], the general solution of
$-f''+qf=\lambda f$ has the form $f(x)=\sum_{j=1}^n
q_j(x)e^{-\kappa_jx}$, where $q_j(x)$ are polynomials. Hence there
exist constants $a_k$, not all zero, such that $\sum_{k=0}^N
a_kf^{(k)}(x)=0$; so by taking the Laplace transform, and introducing
the initial conditions, we can recover
the rational function $\hat f(s)$; see [7,  p.15].
 We recall that any proper rational function
arises as the transfer function of a linear system that has a finite
matrix $A_n$, so $\hat f (\lambda )=C(\lambda I+A_n)^{-1}B$.\par
\rightline{$\square$}\par
\vskip.05in
\noindent {\bf Example (iii).} In Theorem 10.2 (ii), $q$ is a rational
function of $e^{-t/N}$ and under certain conditions gives rise to case
(iii) of Theorem 10.3. In particular, $q(t)=-2{\hbox{sech}}^2 t$ is
algebro-geometric, has period $2\pi i$, is bounded on $\{ z: \vert\Re
z\vert >r\}$ for all $r>0$ and arises from $\tau
(2t)=1+e^{-2t}$. This potential appears in the theory of solitons [4].\par
\vskip.05in
\noindent {\bf 11. Realising linear systems for elliptic potentials}\par
\vskip.05in
\indent Suppose that $q$ is real, smooth and periodic with period
one; introduce Hill's operator
$-{{d^2}\over{dx^2}}+q(x)$ in $L^2({\bf R})$. Then we introduce the Bloch spectrum, which is
$$S_B=\{ \lambda \in {\bf C}: {\hbox{the general
solution of}}\quad -f''+qf=\lambda f\quad {\hbox{has}}\quad f\in 
L^\infty ({\bf R}; {\bf C}) \}.\eqno(11.1)$$
\noindent One can show that when $q$ is an algebro-geometric
potential, $S_B$ has only finitely many gaps; see [1,6, 24]. So
we suppose that 
$$S_B =[\lambda_0, \lambda_1]\cup[\lambda_2,
\lambda_3]\cup\dots\cup [\lambda_{2g}, \infty )\eqno(11.2)$$
\noindent with $g$ gaps. The $\lambda_j$ are the points of the simple
periodic spectrum, such that $-f''+qf=\lambda_jf$ has a unique
solution, up to scalar multiples, that is one or two periodic. Let
$\Phi$ be the $2\times 2$ fundamental solution matrix that satisfies
$${{d}\over{dx}}\Phi (x)=\left[\matrix{ 0&1\cr -\lambda
+q(x)&0\cr}\right] \Phi (x),\qquad \Phi (0)=\left[\matrix{1&0\cr
0&1\cr}\right],\eqno(11.3)$$
\noindent and let $\Delta(\lambda )={\hbox{trace}}\, \Phi (1)$
 be the discriminant of Hill's equation. We can characterize $S_B$ as
$\{ \lambda\in {\bf R}: \Delta(\lambda )^2\leq 4\}$, and its
components are known as
the intervals of stability. \par
\indent The (hyper) elliptic curve ${\cal C}:
y^2=-\prod_{j=0}^{2g}(x-\lambda_j)$ has genus $g$, and we can form the
(hyper) elliptic function field ${\cal E}_g={\bf C}(x)[y]$. We therefore have a situation quite
analogous to (2.5) and the Riemann surface ${\cal E}$ of section 7. In
this section, we consider the case of $g=1$, where ${\cal C}$ is an elliptic curve which is
parametrized by ${\cal P}$.  Hochstadt proved that $g=1$ if
and only if $q(x)=c_1+2{\cal P}(x+c_2)$ where $c_1$ and $c_2$ are
constants; see [17].\par
\indent Starting finite matrices, we formulate a version
of the Gelfand--Levitan equation that is appropriate when $\phi
(x)=Ce^{-xA}B$ is periodic. (The equation (9.12) does not converge when the functions are periodic.) A variant of this was used in
[11] to solve the
matrix nonlinear Schr\"odinger equation. Thus we will realise elliptic tau functions from linear systems.\par
\vskip.05in
\noindent {\bf Definition.} (Periodic linear system $(-A,B,C;E)$) 
Let $A,B,C$ and $E$ be finite square matrices of equal size;
let $\varepsilon=\pm 1$, and suppose that $BC=\varepsilon (AE+EA)$, $BE=EB$, $EA=AE$ and
 $\exp 2\pi A=I$. Define $\phi (x)=Ce^{-xA}B$ to be
the scattering function for $(-A,B,C)$  and then introduce
$$W(x,y)=Ce^{-xA}\bigl(
I-e^{-xA}Ee^{-xA}\bigr)^{-1}e^{-yA}B.\eqno(11.4)$$
\noindent We define the tau function to be 
$$\tau (x)=\exp \Bigl(\int_0^{x} {\hbox{trace}}\,
W(y,y)\, dy\Bigr)\eqno(11.5)$$
\noindent and let $q(x)=-2{{d^2}\over{dx^2}}\log \tau (x)$ be the
potential function.\par
\indent We define $\tau$ indirectly so as to accommodate
the most significant applications. The definition retains the 
spirit of Theorem 9.2, on account of the following result.\par
\vskip.05in
\noindent {\bf Lemma 11.1} {\sl (i) The matrices satisfy the
Gelfand--Levitan equation}
$$-\phi (x+y) +W(x,y)-\varepsilon\int_x^{2\pi} W(x,z)\phi (z+y)\,
dz=W(x,y)E\qquad
(0<x<y<2\pi ), \eqno(11.6)$$
\noindent {\sl and}
$${{d}\over{dx}}\log\det (I-e^{-xA}Ee^{-xA})=\varepsilon {\hbox{trace}}\, W(x,x)
.\eqno(11.7)$$
\indent {\sl (ii) Let ${\bf F}$ be a differential field that contains
all the entries of $e^{-xA}$. Then $\tau (x)$
belongs to a Liouvillian extension of ${\bf F}$, and} $\tau (x+2\pi
)=\kappa \tau (x)$ where $\kappa=\exp\int_0^{2\pi}
{\hbox{trace}}\,W(y,y)\, dy.$
\indent  {\sl (iii) Suppose moreover that $\varepsilon=1$ and $2\pi
\Vert\phi\Vert_\infty <1$. Then}
$${{\partial^2W}\over{\partial x^2}}-{{\partial^2 W}\over{\partial
y^2}}=-2\Bigl({{d}\over{dx}}W(x,x)\Bigr)W(x,y).\eqno(11.8)$$
\vskip.05in
\noindent {\bf Proof.} (i) One can check that 
$$\int_x^{2\pi} e^{-zA}BCe^{-zA}\,
dz=\varepsilon e^{-xA}Ee^{-xA}-\varepsilon E\eqno(11.9)$$
\noindent and it is then a simple matter to verify the integral equation (11.6).\par
\indent By rearranging terms, one checks that
$$\eqalignno{{\hbox{trace}}\,W(x,x)&={\hbox{trace}}\,\bigl(
(I-e^{-xA}Ee^{-xA})^{-1}e^{-xA}BCe^{-xA}\bigr)\cr
&=\varepsilon {{d}\over{dx}}{\hbox{trace}}\log (I-e^{-xA}Ee^{-xA})\cr
&=\varepsilon {{d}\over{dx}}\log\det (I-e^{-xA}Ee^{-xA}).&(11.10)\cr}$$
\indent (ii) By (i), $\tau$ is given by exponential integrals of the entries
of $e^{-xA}$. Note that $W(x,y)$ is periodic in both $x$ and $y$, so
$W(x,x)$ is periodic and hence $\int_0^x{\hbox{trace}}W(y,y)\, dy$
increases by the same amount as $x$ increases through any
interval of length $2\pi$.\par  
\indent (iii)  By repeatedly differentiating (11.6), and using
periodicity, one derives the identity
$${{\partial^2W}\over{\partial x^2}}-{{\partial^2 W}\over{\partial
y^2}}+2\Bigl({{d}\over{dx}}W(x,x)\Bigr)\phi (x+y)+W(x,0)\phi'(y)-{{\partial
W}\over{\partial y}}(x,0)\phi (y)$$
$$-\int_x^{2\pi}\Bigl(
{{\partial^2W}\over{\partial x^2}}-{{\partial^2 W}\over{\partial
y^2}}\Bigr) \phi (z+y)\, dz={{\partial^2W}\over{\partial
x^2}}E-{{\partial^2 W}\over{\partial
y^2}}E\eqno(11.11)$$
\indent Since $ABC-CBA=0$, we obtain 
$$W(x,0)\phi'(y)-{{\partial W}\over{\partial y}}(x,0)\phi
(y)=0,\eqno(11.12)$$
\noindent so (11.11) is a multiple of the original integral equation by
$-2{{d}\over{dx}}W(x,x).$ By the assumptions on 
$\Vert\phi\Vert_\infty$, the solutions are unique, hence the differential equation is
satisfied.\par
\rightline{$\square$}\par
\vskip.05in
\indent By introducing infinite block matrices, we can extend the
scope of Lemma 11.1. Clearly we can replace $\varepsilon$ in (11.6) by a diagonal 
matrix with blocks of $\pm 1$ entries on the diagonal.\par
\indent One can 
interpret the following result as saying that 
Lam\'e's operator $-{{d^2}\over{dx^2}}+2{\cal P}$ has the scattering
function proportional to $\sin x$. Let  $\omega_1$ and 
$\omega_2$ be the periods, so
that $\omega =\omega_2/\omega_1$ has $\Im \omega >0$; then let
$e_1={\cal P}(\omega_1/2)$, $e_2={\cal P}((\omega_1+\omega_2)/2)$ and
$e_3={\cal P}(\omega_2/2)$; then let Jacobi's modulus be
$m^2=(e_2-e_3)/(e_1-e_3)$ and ${\hbox{q}}=e^{i\omega\pi}.$ To be
specific, we choose
$w_1=2\pi$ and $w_2=2\pi i$. Let $A,B$ and $C$ be the infinite block diagonal matrices with
$2\times 2$ diagonal blocks
$$\eqalignno{A&={\hbox{diagonal}}\bigl[J\bigr]_{n=-\infty}^\infty ,
\qquad\qquad  C=A,&(11.13)\cr
E&={\hbox{diagonal}}\bigl[ {\hbox{q}}^{2\vert
n\vert}I_2\bigr]_{n=-\infty}^\infty,\qquad B=2E.&(11.14)\cr}$$
\vskip.05in
\noindent {\bf Proposition 11.2} {\sl (i) The functions $\phi
(x)=Ce^{-xA}B$ and $W(x,y)$ of (11.4) satisfy the Gelfand--Levitan
equation (11.6) and}
$${\hbox{trace}}\,\phi (x)=4{{1+{\hbox{q}}^2}\over{1-{\hbox{q}}^2}}
\sin x;\eqno(11.15)$$
\indent {\sl (ii) The corresponding tau function is entire, belongs to a
Liouvillian extension of the standard elliptic function field
and satisfies}
$$2{\cal P}(x)=-{{d^2}\over{dx^2}}\log\tau (x).\eqno(11.16)$$
\vskip.05in
\noindent {\bf Proof.} (i) The matrices satisfy $EB=BE, AE=EA$ and
$BC=AE+EA$, so Lemma 11.1(i) applies. Note
that the entries of $E$ are summable, so $E$ defines a trace class
operator, hence the trace exists and a simple calculation gives
(11.15).\par
\indent (ii) Observe also that
$$\eqalignno{\det(I-{\hbox{q}}^{2\vert n\vert} e^{-2xA})&
=\det\left[\matrix{1-{\hbox{q}}^{2\vert
n\vert} \cos 2x&-{\hbox{q}}^{2\vert n\vert} \sin 2x\cr
 {\hbox{q}}^{2\vert n\vert}\sin 2x&1-{\hbox{q}}^{2\vert n\vert} \cos 2x\cr}\right]\cr
&=1-2{\hbox{q}}^{2\vert n\vert}\cos 2x+{\hbox{q}}^{4\vert
n\vert},&(11.17)\cr}$$
\noindent so one has 
$$\det (I-e^{-xA}Ee^{-xA})=4\sin^2x\prod_{n=1}^\infty (1-2{\hbox{q}}^{2n}\cos
2x+{\hbox{q}}^{4n})^2;\eqno(11.18)$$
\noindent for comparison, by [25, p 135] the Jacobi elliptic function satisfies
$$\theta_1(x)=2{\hbox{q}}^{1/4}\sin x\prod_{n=1}^\infty \bigl( 1-2{\hbox{q}}^{2n} \cos 2x
+{\hbox{q}}^{4n}\bigr) (1-{\hbox{q}}^{2n}).\eqno(11.19)$$
\noindent So we have an entire function
$$\tau (x)=\det
(I-e^{-xA}Ee^{-xA})={{\theta_1(x)^2}\over{{\hbox{q}}^{1/2}\prod_{n=1}^\infty
(1-{\hbox{q}}^{2n})^2}}.\eqno(11.20)$$
\noindent Moreover, we have [25, p. 132]
$${\cal
P}(x)=-{{d^2}\over{dx^2}}\log\theta_1(x)+e_1+{{d^2}\over{dx^2}}
\log\theta_1(x)\bigr\vert_{x=1/2},\eqno(11.21)$$
\noindent hence we obtain (11.16). Let ${\cal E}$ be the elliptic function field of functions
of rational character on the complex torus ${\bf C}/(\omega_1{\bf
Z}+\omega_2 {\bf
Z})$. Then
${\cal E}={\bf C}({\cal P})[{\cal P}']$, and  by
(11.21) ${\cal
E}$ has a Liouvillian extension ${\cal E}_\theta$ that contains $\theta_1$.

\rightline{$\square$}\par
\vskip.05in
\noindent  {\bf Theorem 11.3} {\sl Let $\tau$ be an elliptic function.
Then there exists a periodic linear
system $(-A,B,C;E)$, where $A,B,C$ and $E$ are infinite block
diagonal matrices with $2\times 2$ blocks, such that} 
$${{d}\over{dx}}\log\tau (x)={\hbox{trace}}\,
W(x,x).\eqno(11.22)$$
\vskip.05in
\noindent {\bf Proof.} Any elliptic function
is of rational character on  ${\bf C}/(\omega_1{\bf
Z}+\omega_2 {\bf Z})$,  and is the ratio of theta functions by [25,
p 105], so
$$\tau (x)=\prod_{j=1}^m {{\theta (x-a_j)}\over{\theta
(x-b_j)}}\eqno(11.23)$$
\noindent where $a_1+\dots +a_m=b_1+\dots +b_m$.\par
\indent First we construct a periodic linear system with $\theta$ as its tau
function. For $n=0$, let
$A_0=J/2$, $E_0=-iJ$, $B_0=iI$ and $C_0=I$, then $(-A_0, B_0, C_0;
E_0)$ is a periodic linear system such that $\det (I-e^{-xA_0}E_0e^{-xA_0})=2i\sin x.$\par
\indent For $n=1, 2, \dots,$ let $A_n=C_n=J$, $E_n={\hbox{q}}^{2n}I$ and
$B_n=2E_n$; then $(-A_n, B_n, C_n;E_n)$ is a periodic linear 
system such that $\det
(I-e^{-xA_n}E_ne^{-xA_n})=1-2{\hbox{q}}^{2n}\cos 2x +{\hbox{q}}^{4n}$.
Hence we can introduce block diagonal matrices 
$A={\hbox{diagonal}}[A_0, A_1, \dots ]$ and $E={\hbox{diagonal}}[E_0,
E_1, \dots ]$, and so on to give a periodic linear system $(-A, B,C;E)$ such that
$$\eqalignno{\det (I-e^{-xA}Ee^{-xA})&=2i\sin x\prod_{n=1}^\infty
(1-2{\hbox{q}}^{2n}\cos 2x +{\hbox{q}}^{4n})\cr
&={{i\theta (x)}\over{{\hbox{q}}^{1/4}\prod_{n=1}^\infty 
(1-{\hbox{q}}^{2n})}}.&(11.24)\cr}$$ 
\indent Next we replace $(-A, B,C;E)$ by the terms $(-A,
e^{a_jA}B, Ce^{a_jA}; e^{a_jA}Ee^{a_jA})$ which give $W_j$ by (11.4);
likewise we introduce 
$(-A,e^{b_jA}B, -Ce^{b_jA}; e^{b_jA}Ee^{b_jA})$ which give $\hat W_j$ by 
(11.4). We then form the block diagonal matrix
$$\oplus_{j=1}^m \Bigl( (-A)\oplus (-A), 
e^{a_jA}B \oplus e^{b_jA}B, Ce^{a_jA}\oplus (-Ce^{b_jA});
 e^{a_jA}Ee^{a_jA}\oplus e^{b_jA}Ee^{b_jA}\Bigr)\eqno(11.25)$$
\noindent which gives the required 
$W(x,y)=\oplus_{j=1}^m W_j(x,y)\oplus \hat W_j(x,y)$ by (11.4), and we verify  
$$\eqalignno{{\hbox{trace}}\,
W(x,x)&=\sum_{j=1}^m\Bigl( {\hbox{trace}}\, W_j(x,x)+{\hbox{trace}}\,
\hat W_j(x,x)\Bigr)\cr
&={{d}\over{dx}}\sum_{j=1}^m \Bigl(\log \theta (x-a_j)-\log\theta
(x-b_j)\Bigr)\cr
&={{d}\over{dx}}\log\tau (x).&(11.26)\cr}$$ 
\noindent  One can
check that $W$ satisfies (11.6) with $\varepsilon$ replaced by a
diagonal matrix with diagonal entries $\pm 1$.\par
\rightline{$\square$}\par
\vskip.05in
\noindent {\bf 12. Linear systems for potentials on hyperelliptic
curves}\par
\vskip.05in
\noindent In this final section, we extend the analysis of section 11 to
(11.2) in the case of
$g>1$. To obtain a model for the Riemann surface of ${\cal C}$, 
we choose a two-sheeted cover of ${\bf C}$ with cuts along $S_B$,
and introduce the canonical homology basis consisting of:\par
\indent $\bullet$ loops $\alpha_j$ that
start from $[\lambda_{2g}, \infty )$, pass along the top sheet to
$[\lambda_{2j-2}, \lambda_{2j-1}]$, then return along the bottom sheet
to the start on  $[\lambda_{2g}, \infty )$;\par
\indent $\bullet$ loops $\beta_j$ that go around the
intervals of stability $[\lambda_{2j-2},
\lambda_{2j-1}]$ that do not intersect with one another, for $j=1,
\dots , g.$ \par
\noindent Then as in [14, p 61], we form the $g\times 2g$ Riemann 
matrix $[I;\Omega]$ from the $g\times g$ matrix blocks
$$I=\Bigl[\int_{\alpha_k}{{x^{j-1}dx}\over{y}}\Bigr]_{j,k=1}^g,
\quad
{\hbox{and}}\quad \Omega=\Bigl[\int_{\beta_k}{{x^{j-1}dx}\over{y}}\Bigr]_{j,k=1}^g.
\eqno(12.1)$$
\noindent Then the corresponding Riemann theta function is
$$\Theta (s\mid \Omega )=\sum_{n\in {\bf Z}^g}\exp\bigl( i\pi \langle
\Omega n,n\rangle +2\pi i \langle s,n\rangle\bigr).\eqno(12.2)$$
\vskip.05in
\noindent {\bf Example.} Suppose that $g=2$ and let 
$$\Omega =\left[\matrix{a&b\cr b&d\cr}\right]\eqno(12.3)$$
\noindent where $\Im a>0, \Im d>0$ and $b\in {\bf Q}$. Then choose
$p\in {\bf N}$ such that $pb\in {\bf Z}$. One can easily check that
$$\Theta (s,t\mid \Omega )=\sum_{r,\mu =0}^{p-1}e^{\pi
(ar^2+2br\mu +d\mu^2)}e^{2\pi rs}e^{2\pi \mu t}\theta (ps+r\mid
p^2a)\theta (pt+\mu \mid p^2d).\eqno(12.4)$$
\vskip.05in
\noindent {\bf Proposition 12.1} {\sl Suppose that $q$ is a
periodic potential with $g$ spectral gaps, as above, and that $\Theta (\,
\,\mid
\Omega )$
is a finite sum of products of Jacobi elliptic functions. Then there 
exist $N<\infty$, $x_j\in {\bf R},$ $\sigma_j\in {\bf C}$ with
$\Im \sigma_j>0$; and block diagonal matrices $(A_j,B_j, C_j; E_j)$ with $2\times
2$ diagonal blocks for $j=1, \dots , N$, such that $\theta (x-x_j\mid \sigma_j)$ is the tau function of
$(-A_j,B_j,C_j;E_j)$ and $q$ belongs to the field ${\bf C}(\theta
(x-x_j\mid \sigma_j); j=1,\dots ,N)$.}\par 
 
\vskip.05in
\noindent {\bf Proof.} We introduce the product of real ovals 
$${\bf T}^g=\Bigl\{ {{1}\over{2}}\Bigl(
\Delta (x_j)+\sqrt{4-\Delta (x_j)^2}\Bigr)_{j=1}^g:\lambda_{2j-1}\leq x_j\leq \lambda_{2j}: j=1, \dots ,
g\Bigr\}\eqno(12.5)$$
\noindent which has dimension $g$. McKean and van Moerbeke [24, p260]
considered the manifold ${\cal M}$ of all the smooth real one-periodic potentials such that
the corresponding Hill's operator has simple spectrum $\{
\lambda_1, \dots , \lambda_{2g}\}$. Using the KdV hierarchy, they 
introduced a $1$ to $1$ differentiable map from ${\cal M}$ onto
 ${\bf T}^g$, where
the tangent vectors on ${\cal M}$ are differential operators. 
Let ${\Lambda}$ be the lattice generated by the columns of $[I; \Omega
]$, and note that ${\bf C}^g/\Lambda$ is the Jacobi variety of ${\cal
C}$. They showed that $q$
extends to an abelian function on ${\bf C}^g$ which is
periodic with respect to $\Lambda$, hence gives a function of rational
character on ${\bf C}^g/\Lambda$. 
The extended function $q$ belongs to ${\cal E}_g$, hence is a theta quotient. 
Moreover, translation on
the potential is equivalent to linear motion on 
${\bf C}^g$ at
constant velocity. \par
\indent 
Thus they solved the inverse spectral problem
explicitly, by showing on [24, p.262] that 
$$q(x)=\sum_{j=0}^{g}\varepsilon_j{{\Theta (X-\omega^*_j/2\mid
\Omega )\Theta
(X-\omega_j^{**}/2\mid\Omega ) }\over{\Theta
(X-\omega_\infty^*/2\mid\Omega )\Theta
(X-\omega_\infty^{**}/2\mid \Omega )}}\eqno(12.6)$$
\noindent where $X=(x_1, \dots , x_{g-1}, ax+b)$ has $a,b,x_1, \dots
,x_{g-1}$ fixed, while $x$ varies, and the constants $\varepsilon_j, 
\omega_j^*, \omega_j^{**}, \omega_\infty^*$ and
$\omega_\infty^{**}$ are notionally computable.\par
\indent By hypothesis, each factor $\Theta (X-\omega^*/2\mid \Omega )$ may be written
a a finite sum of products of functions such as $\theta (ax+c_j\mid
d_j)$, and we can apply Theorem 11.3 to each such factor.\par
\rightline{$\square$}\par 
\vskip.05in
\indent Weierstrass and Poincar\'e developed a systematic reduction
procedure for such elliptic functions of higher genus, so we can
describe the scope of Proposition 12.1. The Siegel upper half-space is 
$${\cal S}_g=\{ \Omega\in M_{g\times g}({\bf C}): \Omega =\Omega^t; \Im
\Omega >0\} .\eqno(12.7)$$
\noindent Let $X$ and $J$ be the $2g\times 2g$ rational block
matrices
$$X=\left[\matrix{ \alpha &\beta \cr \gamma &\delta\cr}\right] ,\qquad
J=\left[\matrix{ 0&-I\cr I&0\cr}\right]\eqno(12.8)$$
\noindent such that $XJX^t=J$; the set of all such $X$ is 
the symplectic group $Sp(2g; {\bf Q}).$ Now $X$ is associated with 
the transformation $\varphi_X$ of ${\cal S}_g$ given by 
$$\varphi_X(\Omega )=(\alpha I+\beta\Omega )^{-1}(\gamma I+\delta \Omega
),\eqno(12.9)$$
\noindent thus $Sp(2g; {\bf Q})$ acts on ${\cal S}_g$. 

\vskip.05in
\noindent {\bf Proposition 12.2} [1] {\sl (i) Suppose that $\Omega$ can be reduced to a
diagonal matrix by the action of the symplectic group. Then
$\Theta (\, \mid \Omega )$ can be expressed as a sum of products of
Jacobian elliptic theta functions.\par
\indent (ii) Condition (i) is equivalent to ${\cal C}$ being a
$N$-sheeted covering of the one-dimensional complex torus for some
$N$.\par
\indent (iii) The orbit of $Sp(2g, {\bf Q})$ that contains $iI$ is dense
in ${\cal S}_g$.}\par
\vskip.05in

\noindent {\bf References}\par
\vskip.05in
\noindent 1. E.D. Belokolos and V.Z. Enolskii, Reduction of abelian
functions and algebraically integrable systems I. Complex analysis and
representation theory. J. Math. Sci. (New York) 106 (2001),
3395--3486.\par 
\noindent 2. M. Bertola, B. Eynard, and J. Harnad, Semiclassical
orthogonal polynomials, matrix models and isomonodromic tau functions,
Comm. Math. Phys. 263 (2006), 401--437.\par 
\noindent 3. G. Blower, Integrable operators and the squares of Hankel
operators, J. Math. Anal. Appl. 340 (2008), 943--953.\par 
\noindent 4. G. Blower, Linear systems and determinantal random point
fields, J. Math. Anal. Appl. 355 (2009), 311--334.\par
\noindent 5. A. Boutet de Monvel, L. Pastur and M. Shcherbina, On the
statistical mechanics approach in the random matrix theory: integrated
density of states, J. Statist. Physics 79 (1995), 585--611.\par
\noindent 6. Y.V. Brezhnev, What does integrability of finite-gap or
soliton potentials mean?, Philos. Trans. R. Soc. Lond. Ser. A
Math. Phys. Eng. Sci. 366 (2008), 923--945.\par
\noindent 7. C.-T. Chen, Linear system theory and design, third edition,
Oxford University Press, 1999.\par
\noindent 8. Y. Chen and A.R. Its, A Riemann--Hilbert approach to the
Akhiezer polynomials, Philos. Trans. R. Soc. Lond. Ser. A
Math. Phys. Eng. Sci. 366 (2008), 973--1003.\par
\noindent 9. Y. Chen and N. Lawrence, A generalization of the
Chebyshev polynomials, J. Phys. A 35 (2002), 4651--4699.\par
\noindent 10. P. Deift, T. Kriecherbauer and K.T. -R. McLaughlin, New
results on the equilibrium measure for logarithmic potentials in the presence of an external
field, J. Approx. Theory 95 (1998), 388--475.\par
\noindent 11. F. Demontis and C. van der Mee, Explicit solutions of the
cubic matrix nonlinear Schr\"odinger equation, Inverse Problems 24 
(2008), 025020.\par    
\noindent 12. L.A. Dikij and I. M. Gelfand, Integrable nonlinear equations
and the Liouville theorem, Funct. Anal. Appl. 13 (1979), 6--15.\par
\noindent 13. H.M. Farkas and I. Kra, Riemann surfaces,
Springer-Verlag, New York, 1980.\par
\noindent 14. A.S. Fokas, A.R. Its, A.A. Kapaev and V. Yu. Novkshenov,
Painlev\'e transcendents: the Riemann--Hilbert approach, American
Mathematical Society, Providence Rhode Island, 2006.\par
\noindent 15. P.L. Forrester and N.S. Witte, Random matrix theory and the
sixth Painlev\'e equation, J. Phys. A 39 (2006), 12211--12233.\par
\noindent 16. R. Fuchs, \"Uber lineare homogone Differentialgleichungen
zweiter Ordnung mit drei in Endlichen gelegenen wesentlich
 singul\"aren Stellen, Math. Ann. 63 (1907),
301--321.\par 
\noindent 17. F. Gesztesy and R. Weikard, Elliptic algebro-geometric
solutions of the KdV and AKNS hierarchies--an analytic approach,
 Bull. Amer. Math. Soc. (N.S.) 35 (1998), 271--317.\par 
\noindent 18. D. Guzzetti, The elliptic representation of the general
Painlev\'e VI equation, Comm. Pure Appl. Math. 55 (2002), 1280--1363.\par
\noindent 19. N.J. Hitchin, Riemann surfaces and integrable systems, pp
11--52 in N.J. Hitchin, G.B. Segal and R.S. Ward, Integrable systems:
twistors, loop groups and Riemann surfaces, Oxford Science
Publications, 1999.\par
\noindent 20. M. Jimbo, T. Miwa and K. Ueno, Monodromy preserving
deformation of linear ordinary differential equations with rational
coefficients, I. General theory and $\tau$ function, Phys. D 2 (1981),
306--352.\par
\noindent 21. K. Johansson, On fluctuations of eigenvalues of random
Hermitian matrices, Duke Math. J. 91 (1998), 151--204.\par
\noindent 22. A.P. Magnus, Painlev\'e
-type differential equations for the recurrence coefficients of
semi-classical orthogonal polynomials, J.Comp. Appl. Math. 57 (1995),
215--237.\par
\noindent 23. R. Maier, Lam\'e polynomials, hyperelliptic reduction and 
Lam\'e band structure, 
 Philos. Trans. R. Soc. Lond. Ser. A
Math. Phys. Eng. Sci. 366 (2008), 1115--1153.\par
\noindent 24. H.P. McKean and P. van Moerbeke, The spectrum of Hill's
equation, Invent. Math. 30 (1975), 217--274.\par
\noindent 25. H. McKean and V. Moll, Elliptic curves: Function theory,
geometry, arithmetic, Cambridge University Press, 1997.\par
\noindent 26. M. L. Mehta, Random matrices, second edition (Academic
Press, San Diego 1991).\par
\noindent 27. K. Okamoto, On the $\tau$-function of the Painlev\'e
equations, Physica D 2 (1981), 525--535.\par
\noindent 28. L.A. Pastur, Spectral and probabilistic aspects of matrix
models, pp 205--242, in Algebraic and Geometric Methods in
Mathematical Physics, edrs. A. Boutet de Monvel and V.A. Marchenko,
Kluwer Acad. Publishers, 1996.\par
\noindent 29. E. B. Saff and V. Totik, Logarithmic potentials with
external fields, Springer, Berlin 1997.\par
\noindent 30. M.F. Singer, Introduction to the Galois theory of linear differential
equations, pp. 1--83, in {\sl Algebraic theory of Differential
Equations}, Edrs M.A.H. McCallum and A.V. Mikhailov, London
mathematical Society Lecture Notes, Cambridge, 2009.\par
\noindent 31. C.A. Tracy and H. Widom, Fredholm determinants,
differential equations and matrix models, Comm. Math. Phys. 163
(1994), 33--72.\par
\vfill
\eject
\end